\documentclass[12pt]{amsart}
\usepackage{amssymb,amsmath,amsfonts,latexsym}
\usepackage{bm, enumerate}
\usepackage{hyperref, xcolor}

\newcommand{\newsection}[1]{\setcounter{equation}{0} \section{#1}}
\newcommand{\bea}{\begin{eqnarray}}
\newcommand{\eea}{\end{eqnarray}}

\newcommand{\clb}{\mathcal{B}}

\newcommand{\cld}{\mathcal{D}}

\newcommand{\clh}{\mathcal{H}}

\newcommand{\clm}{\mathcal{M}}

\newcommand{\D}{\mathbb{D}}
\newcommand{\N}{\mathbb{N}}
\newcommand{\C}{\mathbb{C}}

\def\textmatrix#1&#2\\#3&#4\\{\bigl({#1 \atop #3}\ {#2 \atop #4}\bigr)}
\def\dispmatrix#1&#2\\#3&#4\\{\left({#1 \atop #3}\ {#2 \atop #4}\right)}
\newcommand{\be}{\begin{equation}}
\newcommand{\ee}{\end{equation}}
\newcommand{\ben}{\begin{eqnarray*}}
\newcommand{\een}{\end{eqnarray*}}

\newcommand{\NI}{\noindent}

\newcommand{\bi}{\begin{itemize}}
\newcommand{\ei}{\end{itemize}}

\newcommand\la{\langle}
\newcommand\ra{\rangle}

\newtheorem{Theorem}{\sc Theorem}[section]
\newtheorem{Lemma}[Theorem]{\sc Lemma}
\newtheorem{Proposition}[Theorem]{\sc Proposition}
\newtheorem{Corollary}[Theorem]{\sc Corollary}
\newtheorem{Definition}[Theorem]{\sc Definition}
\newtheorem{Question}{\sc Question}
\newtheorem{ass}[Theorem]{\sc Assumption}

\theoremstyle{definition}
\newtheorem{Example}[Theorem]{\sc Example}
\newtheorem{Remark}[Theorem]{\sc Remark}

\newtheorem{Note}[Theorem]{\sc Note}

\newcommand{\bt}{\begin{Theorem}}
\def\beginlem{\begin{Lemma}}
\def\beginprop{\begin{Proposition}}
\def\begincor{\begin{Corollary}}
\def\begindef{\begin{Definition}}
\def\beginexamp{\begin{Example}}
\def\beginrem{\begin{Remark}}
\def\beginq{\begin{Question}}
\def\beginass{\begin{ass}}
\def\beginnote{\begin{Note}}
\newcommand{\et}{\end{Theorem}}
\def\endlem{\end{Lemma}}
\def\endprop{\end{Proposition}}
\def\endcor{\end{Corollary}}
\def\enddef{\end{Definition}}
\def\endexamp{\end{Example}}
\def\endrem{\end{Remark}}
\def\endq{\end{Question}}
\def\endass{\end{ass}}
\def\endnote{\end{Note}}

\numberwithin{equation}{section}

\begin{document}

\title{Completely non-unitary contractions and analyticity }

\author[Das]{Susmita Das}
\address{Indian Institute of Science, Department of Mathematics, Bangalore, 560012,
India}
\email{susmitadas@iisc.ac.in, susmita.das.puremath@gmail.com}

\date{\today}


\begin{abstract}
We study completely non-unitary contractions $T$ with finite dimensional defect spaces $\cld_T$ and $\cld_{T^*}$. We present a complete classification of all such contractions $T$ that satisfy a generalized property of Hardy shift operator: $\cld_T\subseteq\cld_{T^*}$. We show that $T$ is analytic if and only if it has no non-zero eigenvalue. Furthermore, we characterize the hyponormality of all those $T$.
\end{abstract}

\maketitle

\newsection{Introduction}\label{sec: intro}

Let $\clh$ be a separable, infinite-dimensional, complex Hilbert space, and $\clb(\clh)$ denote the algebra of all bounded linear operators on $\clh$. An operator $T \in \clb(\clh)$ is called a contraction if $\|T\|\leq 1$. To each contraction $T$, one can associate two operators, $D_T = (I-T^*T)^{1/2}$
and $D_{T^*} = (I-TT^*)^{1/2}$, called the defect operators (\cite{NF}) -- which are trivial when $T$ is unitary. The spaces $\cld_T = \overline{\cld_T \clh}$ and $\cld_{T^*}=\overline{D_{T^*}\clh}$ are known as the defect spaces and their corresponding dimensions are called the defect indices of $T$. In this paper, we are interested in studying a class of operators having the property $\cld_T\subseteq\cld_{T^*}$. Hardy shift operator, which we denote by $S$ throughout, trivially possesses this property, as $\cld_S=0$. Another instance, where this property arises is with hyponormal contractions, which we discuss below (also to motivate the problems we solve in this paper).
Recall, $T\in \clb(\clh)$ is hyponormal (\cite{Martin-Putinar}), if $T^*T\geq TT^*$. For such $T$, the range, Ran~$T$ is always contained in $\text{Ran }T^*$ by Douglas Lemma (\cite{Douglas}), which asserts the equivalence of the following statements about any two operators $A, B\in \clb(\clh)$ as:
\textsf{\begin{enumerate}
\item $\text{Ran }A \subseteq\text{Ran }B $.
\item $AA^*\leq \lambda^2 BB^* \text{ for some }\lambda\geq 0$.
\item \text{There exists a bounded operator } $C$ on $\clh$ such that $A=BC$.
\end{enumerate}
}
The properties of hyponormal $T$ and their spectral behaviors have been extensively studied in the literature (e.g. \cite{Curto 1}, \cite{Curto 2}, \cite{Curto}, \cite{DB}, \cite{Exner}, \cite{Xia}). However, a complete characterization of hyponormality for a general bounded linear operator is still unknown. Now if a hyponormal $T$ is also a contraction, one is helped by Douglas lemma. Namely, for such $T$, the relation $T^*T\geq TT^*$ and the existence of defect operators together yield $I-T^*T\leq I-TT^*$, which is exactly the second condition of the Douglas lemma for $A=D_T$, $B=D_{T^*}$ and $\lambda=1$. Consequently, hyponormal contractions $T$ will always satisfy $\cld_T\subseteq\cld_{T^*}$.

Beginning with the classes of operators $T$ of our interest, a contraction $T$ is called completely non-unitary (c.n.u. for short), if there is no nontrivial $T$ reducing subspace $\clm$ of $\clh$ such that $T|_\clm$ is unitary. It is well known (\cite{NF}), that any nontrivial bounded contraction can be uniquely decomposed as the direct sum of a unitary and a completely non-unitary contraction. Throughout the paper, we will focus on c.n.u. contractions, which suffice, as the unitary part of a contraction have no role in the computation of the associated defect operators and the corresponding indices. The Unilateral shift $S$ on the Hardy space $H^2(\D)$ is one of the well-studied examples of c.n.u. contractions (\cite{AS}). Now as $S$ is an isometry, $D_S=0$ and so $\cld_S\subseteq\cld_{S^*}$, and moreover $D_{S^*}$ is of finite dimension. Observe that, $S$ is also hyponormal. These properties raise the natural question: for a c.n.u contraction $T$ with finite indices and satisfying $\cld_T\subseteq\cld_{T^*}$, can we recover hyponormality? The answer is `No'. In Section \ref{sec:Characterization}, we give a counter-example, which shows that these operators are not hyponormal in general and this leads us to several problems that we solve in this paper, namely, first understanding the class of operators $T$ with the aforementioned properties and then determining the hyponormal operators among them, and studying their further properties.

 In Theorem \ref{thm: char}, we first characterize a c.n.u contraction $T$ with finite indices and having the property $\cld_T\subseteq\cld_{T^*}$. This characterization result establishes a deep connection between the theory of contraction and the theory of finite rank perturbations of unilateral shift (\cite{Nakamura 1}, \cite{Nakamura}). Based on this, we completely classify the hyponormality of such contractions in terms of certain finite rank operators in Theorem \ref{Hyponormal}. As an application, we show in Theorem~\ref{D_T=1} that, if in addition, $T$ is a c.n.u. with $\dim\cld_T=1$, $T$ will always be hyponormal. Now, let us turn to the analyticity part.

A bounded linear operator $T$ on a Hilbert space $\clh$ is called analytic if $\bigcap_{n\geq 1}T^n \clh=\{0\}$. A wide class of analytic operators includes all nilpotent operators, and the multiplication operators $M_z$ on reproducing kernel Hilbert spaces of complex analytic functions (\cite{Aronszajn}). The property of analyticity plays a crucial role in (\cite{SS}) to establish the unitary equivalence of a left invertible operator and the multiplication operator $M_z$ on some reproducing kernel Hilbert space. Here is a curious note: every analytic contraction is a c.n.u. However, a c.n.u. contraction $T$ with finite indices neither satisfies $\cld_T\subseteq\cld_{T^*}$, nor it has to be analytic in general (which motivates the next problem we solve). To quickly verify this, consider the bounded operator $T$ on $H^2(\D)$ defined by, $$T(1)=T(z)=\frac{1}{\sqrt{2}},\quad \text{and }\ \  T(z^n)=z^{n+1}\  \text{ for all } n\geq 2.$$ Then, an easy computation shows that $T$ is a c.n.u. and $\cld_T=\mathbb{C}(1-z)$, $\cld_{T^*}=\C z\oplus \C z^2$, and hence $\cld_T\nsubseteq\cld_{T^*}$. Next, as $1$ is an eigenvector of $T$, it is not analytic. Under these circumstances, one can seek to classify all the analytic c.n.u. contractions. We provide an answer for the analyticity in the light of our characterization result, Theorem~\ref{thm: char}.
\medskip

\noindent
The rest of the paper is organized as follows:
In Section \ref{sec: example}, we consider an example of an operator $T$ which is not an isometry, but a c.n.u contraction with $\cld_T \subseteq \cld_{T^*}$ and $\dim\cld_{T^*}<\infty$. We eventually show that $T$ is both hyponormal and analytic. In Section \ref{sec:Characterization} we provide a complete characterization of c.n.u. $T$ with finite indices and $\cld_T\subseteq\cld_{T^*}$ (Theorem \ref{thm: char}). We show that these contractions are certain finite rank perturbations of unilateral shifts of finite multiplicities. Using this structure theorem, we classify all the c.n.u. contractions that are hyponormal in Theorem \ref{Hyponormal}. In the last section \ref{sec: Analyticity}, we prove that such c.n.u.  $T$ is analytic, if and only if it has no non-zero eigenvalue.

\newsection{An Example}\label{sec: example}

We begin by recalling once again, the unilateral shift $S$ on the Hardy space $H^2(\D)$ has finite-dimensional defect spaces; explicitly, $\cld_S =\{0\}$ (since $S$ is an isometry) and $\cld_{S^*}=\C1$ so that, $\cld_ S\subseteq \cld_{S^*}$ holds trivially. Furthermore, $S$ is a c.n.u contraction and is also analytic i.e., $\bigcap_{n\geq 1}S^n H^2(\D)=\{0\}$.
Let us now look at an example wherein the defect spaces are non-trivial and finite-dimensional (and with $\cld_T\subseteq \cld_{T^*}$):
\begin{Example}
First fix two non-zero complex numbers $a_1, a_2$ with $|a_1|^2+ |a_2|^2 <1$. Consider the rank one perturbation $T$ of the square of unilateral shift operator on $H^2(\D)$ defined by $$T(1)= a_1z +a_2 z^2\quad  \text{and}\quad T(z^n)= z^{n+2}\ \text{ for } n\geq 1.$$
With respect to the orthonormal basis $\{z^n\}_{n\geq 0}$ on $H^2(\D)$, $T$ and $T^*$ have the matrix representations:
\begin{equation}\label{eqn: T matrix}
[T] = \begin{bmatrix}
0 & 0 & 0 & 0 & \cdots
\\
a_1& 0 & 0 & 0 & \cdots
\\
a_2 & 0 & 0 & 0 & \cdots
\\
0 & 1 & 0 & 0 & \cdots
\\
0 & 0 & 1 & 0 & \cdots
\\
0 & 0 & 0 & 1 & \cdots
\\
\vdots & \vdots & \vdots& \vdots &\ddots
\end{bmatrix}, \qquad
[T^*] = \begin{bmatrix}
0 & \bar{a_1} & \bar{a_2} & 0 & 0 & 0 & \cdots
\\
0 & 0 & 0 & 1 & 0 & 0 &\cdots
\\
0 & 0 & 0 & 0 & 1 & 0 &\cdots
\\
0 & 0 & 0 & 0 & 0 & 1 & \cdots
\\
\vdots & \vdots & \vdots &\vdots &\vdots &\vdots &\ddots
\end{bmatrix}.
\end{equation}

A simple computation shows (again with respect to the same basis),

\begin{equation}\label{eqn: D_T matrix}
[I-T^*T] = \begin{bmatrix}
1-(|a_1|^2+|a_2|^2)& 0 & 0 & \cdots
\\
0 & 0 & 0  & \cdots
\\
0 & 0 & 0 & \cdots
\\
0 & 0 & 0 & \cdots
\\
0 & 0 & 0 & \cdots
\\
\vdots & \vdots & \vdots&\ddots
\end{bmatrix},\
[I-TT^*] = \begin{bmatrix}
1 & 0 & 0 & 0 & 0 & \cdots
\\
0 & 1-|a_1|^2 & -a_1\bar{a_2} & 0 & 0 & \cdots
\\
0 & -\bar{a_1}a_2 & 1-|a_2|^2 & 0 & 0 & \cdots
\\
0 & 0 & 0 & 0 & 0 & \cdots
\\
0 & 0 & 0 & 0 & 0 & \cdots
\\
\vdots & \vdots & \vdots & \vdots &\vdots &\ddots
\end{bmatrix}.
\end{equation}

Since, both $1-(|a_1|^2+|a_2|^2)$ and $(1-|a_1|^2)$ are strictly positive, and $$\det\begin{pmatrix} 1 & 0 & 0\\ 0 & 1-|a_1|^2 & -a_1\bar{a_2} \\0 & -\bar{a_1}a_2 & 1-|a_2|^2 & \end{pmatrix}=1-(|a_1|^2+|a_2|^2)>0,$$ it follows that, $T^*T\leq I$ and $TT^*\leq I$. Hence $T$ is a contraction, $\cld_T= \mathbb{C}1$ and $\cld_{T^*}=\overline{span}\{ 1, z, z^2\}$; in particular, $\cld_T\subseteq\cld_{T^*}$.

Now, showing the analyticity of $T$. Indeed, $T(1)=a_1z+a_2 z^2= z(a_1+a_2z)$ and $T(z^n)=z^{n+2}$ for all $n\geq 1$ imply, $T(H^2(\D))\subseteq zH^2(\D)$ and further $T^n(H^2(\D))\subseteq z^{n+2}(H^2(\D))$ for $n\geq 2$. Hence $\bigcap_{n\geq 1}T^n(H^2(\D))\subseteq\bigcap_{n\geq 1}z^n H^2(\D)=\{0\}$.

\smallskip

Every analytic contraction $T$ on a Hilbert space $\clh$ is completely nonunitary. Indeed, if $M\subseteq\clh$ is such that $T|_M$ is unitary, then $M= TM$, and hence $M \subseteq \bigcap_{n \geq 1} T^n\clh =\{0\}$, forcing $M=0$. Since $T$ in our context is analytic, it is also c.n.u.

As we explain, it is also easy to see that $T$ is hyponormal. First, an easy computation shows that the matrix representation of $T^*T-TT^*$ with respect to the basis $\{ z^n\}_{n\geq 0}$, is given by

\begin{equation}\label{eqn: hyponormal matrix}
[T^*T-TT^*] = \begin{bmatrix}
|a_1|^2+|a_2|^2 & 0 & 0 & 0 & 0 & \cdots
\\
0 & 1-|a_1|^2 & -a_1\bar{a_2} & 0 & 0 & \cdots
\\
0 & -\bar{a_1}a_2 & 1-|a_2|^2 & 0 & 0 & \cdots
\\
0 & 0 & 0 & 0 & 0 & \cdots
\\
0 & 0 & 0 & 0 & 0 & \cdots
\\
\vdots & \vdots & \vdots & \vdots&\vdots &\ddots
\end{bmatrix}.
\end{equation}

Since $(|a_1|^2+|a_2|^2)$, $(1-|a_1|^2)$ are strictly positive and $$\det \begin{pmatrix} |a_1|^2+|a_2|^2 & 0 & 0\\ 0 & 1-|a_2|^2 & -a_1\bar{a_2}\\0 & -\bar{a_1}a_2 & 1-|a_2|^2\end{pmatrix}=(|a_1|^2+|a_2|^2)(1-(|a_1|^2+|a_2|^2))>0,$$ $T^*T-TT^*\geq 0$ and hence $T$ is hyponormal.
\end{Example}
For any hyponormal operator $T$, $\cld_T\subseteq\cld_{T^*}$ (as explained earlier); we will soon see that the converse is not always true. In the following section, we aim to obtain a characterization for c.n.u. contractions $T$ with defect spaces having arbitrary finite dimensions.

\newsection{Characterization}\label{sec:Characterization}

We begin with setting up some notations. Corresponding to an orthonormal basis $\{e_m\}_{m\geq 1}$ on a Hilbert space $\clh$, and $n\in \mathbb{N}$, let $P_n$ denote the orthogonal projection onto the span of $\{e_1, e_2, \ldots, e_n\}$.
For $k \geq 1$, $S_k$ denotes the unilateral shift of multiplicity $k$ on $\clh$ i.e.,
\begin{equation}\label{S_k}
S_k(e_m)=e_{m+k},\quad m \geq 1.
\end{equation}
With these notations (and the basis), we define a family of finite rank operators $F$, $F_r$ ($r\in \mathbb{N}$):
We fix some scalars $x^i_{i+k}\in\C$,  $i\in \{1, 2,\cdots, n\}$ and $k\geq 1$, and define
\begin{eqnarray}
 F(e_i)&=&
\begin{cases}
x_1^ie_1+x_2^ie_2+\cdots+(x_{i+k}^i-1)e_{i+k}+\cdots+x_{n+k}^ie_{n+k} & \mbox{if } i\leq n, \\
0 & \mbox{if } i\geq n+1.
\end{cases}\label{F}\\
 F_1&=& F+ S_kP_n, \label{F_1'}\\
 F_r& =& F_1^r+S_k(I-P_n)F_1^{r-1}+\cdots + S_k^{r-1}(I-P_n)F_1 \quad \forall\  r\geq 2.\label{F_r}
\end{eqnarray}
In this setting, we will prove the following proposition which is crucial in our paper:
\begin{Proposition}\label{prop: char}
Let $n \geq 0$ and $ k \in \mathbb{N}$ be arbitrary but fixed. Let $\clh$, $\{e_m\}_{m\geq 1}$, $S_k$, $F$ and $F_r$ ($r\in \mathbb{N}$) be as above \eqref{S_k}--\eqref{F_r}. Set $T=S_k+F$, and assume that $T$ is a c.n.u contraction with finite indices such that $\cld_T\subseteq\cld_{T^*}$, $\dim\cld_{T}=n$ and $\dim(\cld_{T^*}\ominus\cld_T)=k$. Then:
\begin{enumerate}
\item $\text{rank} (P_n-F_1^*F_1)=n$.\smallskip
\item $F_1(I-P_n)=(I-P_{n+k})F_1=0$.\smallskip
\item $\lambda\|F_1^*x\|^2-\|F_1x\|^2\leq \lambda \|P_{n+k} x\|^2-\|P_n x\|^2, \quad \forall x\in \clh$ and for some $\lambda\geq 0$.\smallskip
\item $\|F_rx\|\leq \|P_n x\|$ \text{and} $\|F_r^*x\|\leq \|P_{n+kr}x\|$  \text{hold for all} $r \in \mathbb{N}$ and $x\in \clh$,
with zero as the only common solution to the corresponding equalities.
\end{enumerate}
\end{Proposition}
\begin{proof}
Begin with $T=S_k+F$ as in the statement; recall the definitions of $S_k,\ F,\ F_r$. With respect to the basis $\{e_m\}_{m\geq 1}$ as above, $T$ can be represented by the following matrix:
\begin{equation}\label{eqn: T new matrix}
[T] = \begin{bmatrix}
x^1_1 & x^2_1 &\cdots &x^n_1 & 0 & 0 & 0 & \cdots
\\
x^1_2 & x^2_2 &\cdots & x^n_2 & 0 & 0 & 0 & \cdots
\\
\vdots & \vdots & \vdots & \vdots & \vdots &\vdots &\vdots & \vdots
\\
x^1_k & x^2_k &\cdots & x^n_k & 0 & 0 & 0 & \cdots
\\
1+(x^1_{k+1}-1) & x^2_{k+2} &\cdots & x^n_{k+1} & 0 & 0 & 0 & \cdots
\\
\vdots & \vdots & \vdots &\vdots &\vdots &\vdots & \vdots &\vdots
\\
x^1_{n+k} & x^2_{n+k} &\cdots & 1+(x^n_{n+k}-1) & 0 & 0 & 0 & \cdots
\\
0 & 0 &\cdots & 0 & 1 & 0 & 0 & \cdots
\\
0 & 0 &\cdots & 0 & 0 & 1 & 0 & \cdots
\\
\vdots & \vdots & \vdots &\vdots &\vdots &\vdots & \vdots &\ddots
\end{bmatrix}.
\end{equation}
As $F=F_1-S_kP_n$, we have by equation \eqref{F}:
\[F_1(e_i)=
\begin{cases}
x_1^ie_1+x_2^ie_2+\cdots+x_{i+k}^ie_{i+k}+\cdots+x_{n+k}^ie_{n+k} & \mbox{if } i\leq n, \\
0, & \mbox{if } i\geq n+1.
\end{cases}
\]
This immediately implies the below equation, which proves assertion $(2)$:
\begin{equation}\label{eqn: F_1(I-P_n)=(I-P_{n+k})F_1=0}
F_1(I-P_n)=(I-P_{n+k})F_1=0.
\end{equation}
\smallskip
Again by $F=F_1-S_kP_n$, we can write $T=S_k+F$ as
\begin{equation}\label{T= F_1+ S_k(I-P_n)}
T= F_1+ S_k(I-P_n).
\end{equation}\smallskip
Claim: For $r\geq 2$,
\begin{equation}\label{T^r}
T^r= F_1^r+ S_k(I-P_n)F_1^{r-1}+\cdots +S_k^{r-1}(I-P_n)F_1+S_k^r(I-P_n).
\end{equation}
\smallskip
We prove the claim by induction on $r\geq 2$. The proof when $r=2$, follows by the computation-
\begin{equation*}
\begin{split}
T^2 &= [F_1+S_k(I-P_n)][F_1+ S_k(I-P_n)] \\
&= F_1^2 + F_1 S_k(I-P_n)+S_k(I-P_n)F_1+S_k(I-P_n)S_k(I-P_n) \\
&= F_1^2 +0 +S_k(I-P_n)F_1 + S_k^2(I-P_n) \\
&= F_1^2 + S_k(I-P_n)F_1 + S_k^2(I-P_n).
\end{split}
\end{equation*}
For the induction step, assume the claim to be true for some $r>2$. Then, we are done by-
\begin{equation*}
\begin{split}
T^{r+1} &= T^r T  \\
&= [ F_1^r+ S_k(I-P_n)F_1^{r-1}+\cdots +S_k^{r-1}(I-P_n)F_1+S_k^r(I-P_n)][F_1+S_k(I-P_n)] \\
&= F_1^{r+1}+ S_k(I-P_n)F_1^r+\cdots+ S_k^{r-1}(I-P_n)F_1^2+ S_k ^r(I-P_n)F_1+ S_k^r(I-P_n)S_k(I-P_n) \\
&= F_1^{r+1}+ S_k(I-P_n)F_1^r+ \cdots +S_k^r(I-P_n)F_1 + S_k^{r+1}(I-P_n). \\
\end{split}
\end{equation*}
By the above claim \eqref{T^r} and definition (\ref{F_r}), we deduce that,
\begin{equation}\label{T^r short}
T^r= F_r+ S_k^r(I-P_n) , r \geq 2.
\end{equation}
\smallskip
We now compute the defect operators $\cld_T$ and $\cld_{T^*}$ of $T$. Note that,
\begin{equation*}
\begin{split}
T^*T &= \Big(F_1^*+(I-P_n)S_k^*\Big)\Big(F_1+S_k(I-P_n)\Big) \\
&= F_1^*F_1+ F_1^*S_k(I-P_n)+(I-P_n)S_k^*F_1+(I-P_n)S_k^*S_k(I-P_n) \\
& = F_1^*F_1 + 0 + 0 + (I-P_n),\\ \text{and hence,}
\end{split}
\end{equation*}
\begin{equation}\label{D_T=D_F_1}
I-T^*T= P_n-F_1^*F_1.
\end{equation}
Since $\dim \cld_T=n$, the below equation finishes the proof of the first assertion. $$\text{Rank}(P_n-F_1^*F_1)=n.$$
Again,
\begin{equation*}
\begin{split}
TT^* &=\Big(F_1+S_k(I-P_n)\Big)\Big(F_1^*+(I-P_n)S_k^*\Big) \\
& F_1F_1^* + F_1(I-P_n)S_k^*+S_k(I-P_n)F_1^*+S_k(I-P_n)S_k^*.
\end{split}
\end{equation*}
Now equation \eqref{eqn: F_1(I-P_n)=(I-P_{n+k})F_1=0} implies $F_1(I-P_n)=0$ and so $(I-P_n)F_1^*=0$. Therefore,
\begin{equation}\label{D_{T^*}}
TT^* = F_1F_1^* + S_k(I-P_n)S_k^*.
\end{equation}

We next show that, $$I-S_k(I-P_n)S_k^*=P_{n+k}.$$
Indeed for $1\leq m <k$, $S_k^*(e_m)=0$ and hence in this case, $\Big(I-S_k(I-P_n)S_k^*\Big)e_m= e_m$. For $k\leq m\leq n+k$, we have $S_k^* e_m=e_{m-k}\in \{e_1, e_2, \ldots, e_n\}$ where $e_0=0$. Therefore, $(I-P_n)S_k^*e_m=0$ and $\Big(I-S_k(I-P_n)S_k^*\Big)e_m = e_m$ for $k\leq m\leq n+k$.
Now, for $m \geq n+k+1$, $\Big(I-S_k(I-P_n)S_k^*\Big)e_m= e_m- S_k(I-P_n)e_{m-k}=e_m-S_k e_{m-k}=e_m-e_m=0$. Therefore,
\[\Big(I-S_k(I-P_n)S_k^*\Big)e_m=
\begin{cases}
e_m & \mbox{if } 1\leq m \leq n+k, \\
0 & \mbox{if }  m\geq n+k+1,
\end{cases}
\]
allowing us to write-
\begin{equation}\label{P_{n+k}}
\Big(I-S_k(I-P_n)S_k^*\Big)e_m= P_{n+k}.
\end{equation}
Equations (\ref{D_{T^*}}) and (\ref{P_{n+k}}) together imply
\begin{equation}\label{D_T*=D_F_1^*}
I-TT^*= P_{n+k}-F_1F_1^*.
\end{equation}
For $r\geq 2$, by equation \eqref{T^r short} we have
\begin{equation}\label{T^*r T^r}
\begin{split}
T^{*r} T^r &= \Big(F_r^*+(I-P_n)S_k^{*r}\Big)\Big(F_r+S_k^r(I-P_n)\Big) \\
&= F_r^*F_r+F_r^*S_k^r(I-P_n)+(I-P_n)S_k^{*r}F_r+(I-P_n)S_k^{*r}S_k^r(I-P_n).
\end{split}
\end{equation}
And by definition \eqref{F_r},
\begin{equation*}
F_r^* = F_1^{*r}+F_1^{*(r-1)}(I-P_n)S_k^*+ \cdots +F_1^*(I-P_n)S_k^{*(r-1)}.
\end{equation*}
Hence,
\begin{equation}\label{F_rS_k}
F_r^*S_k^r(I-P_n)=F_1^{*r}S_k^r(I-P_n)+F_1^{*(r-1)}(I-P_n)S_k^{(r-1)}(I-P_n)+ \cdots + F_1^*(I-P_n)S_k(I-P_n).
\end{equation}
Since $F_1^*(e_{n+k+l})=0$ for all $l\geq 1$, the above equation implies $F_r^*S_k^r(I-P_n)=0$, and hence $(I-P_n)S_k^{*r}F_r=0$. This, together with (\ref{T^*r T^r}) yields,
\begin{equation}\label{D_{T^r}}
I-T^{*r}T^r=P_n-F_r^*F_r,\quad \forall\ r\geq 2.
\end{equation}
Again, for $r \geq 2$,
\begin{equation}\label{T^r T^*r}
\begin{split}
T^rT^{*r} &= \Big(F_r+ S_k^r(I-P_n)\Big)\Big(F_r^*+(I-P_n)S_k^{*r}\Big) \\
&= F_rF_r^*+ F_r(I-P_n)S_k^{*r}+S_k^r(I-P_n)F_r^*+ S_k^r(I-P_n)S_k^{*r}.
\end{split}
\end{equation}
Since $F_1(e_m)=0$ for all $m \geq (n+1)$,  it follows by the definition of $F_r$ that, $F_r(I-P_n)S_k^{*r}=0$. Hence, $S_k^r(I-P_n)F_r^*=0$  for  $r \geq 2$.
Again for $l\geq 1$ and $r\geq 2$, $S_k^r(I-P_n)S_k^{*r}(e_{n+kr+l})=S_k^r(e_{n+l})= e_{n+kr+l}$ and for $ 1\leq m \leq (n+kr)$, $S_k^r(I-P_n)S_k^{*r}(e_m)=0$ and hence,
\[S_k^r(I-P_n)S_k^{*r}(e_m)=
\begin{cases}
e_m & \mbox{if } m\geq n+kr+1, \\
0, & \mbox{if }   m\leq n+kr,
\end{cases}\qquad \forall\ r\geq 2.
\]
\begin{equation}\label{P_{n+kr}}
\text{So,}\quad S_k^r(I-P_n)S_k^{*r}=I- P_{n+kr}, \quad \forall\ r\geq 2.
\end{equation}
Therefore it follows by equations \eqref{T^r T^*r} and \eqref{P_{n+kr}},
\begin{equation}\label{D_{T^*r}}
I- T^rT^{*r}= P_{n+kr}-F_rF_r^*, \quad r\geq 2.
\end{equation}
Combining the equations \eqref{D_T=D_F_1}, \eqref{D_{T^r}}, \eqref{D_T*=D_F_1^*} and \eqref{D_{T^*r}}, we have for all $r\geq 1$,
\begin{eqnarray}
D_{T^r}^2&=& (I-T^{*r}T^r)\ =\ (P_n-F_r^*F_r),\ \ \text{and}\label{D_Tr}\\ D_{T^{*r}}^2 &=& (I-T^rT^{*r})\ =\ (P_{n+kr}-F_rF_r^*).\label{D_T*r}
\end{eqnarray}
The implications below \eqref{ker D_T^r} show for every $r\geq 1$,
\begin{equation}\label{ker D_T^r}
\ker D_{T^r}=\{x \in \clh: \|F_r x\|=\|P_n x\|\},
\end{equation}
\begin{align*}
\begin{aligned}
 D_{T^r}x =0
\iff (P_n-F_r^*F_r)^{1/2} x =0
\iff & \la(P_n-F_r^*F_r)^{1/2} x, (P_n-F_r^*F_r)^{1/2} x\ra =0
\\
\iff & \la(P_n-F_r^*F_r) x, x\ra =0
\\
\iff & \la F_r x, F_r x\ra =\la P_n x, x\ra
\\
\iff & \| F_r x\| = \|P_n x\|
\end{aligned}
\end{align*}
Similarly, it can be shown that for $r \geq 1$,
\begin{equation}\label{ker D_T^*r}
\ker D_{T^{*r}}=\{x \in \clh: \|F_r^* x\|=\|P_{n+kr}x\|\}.
\end{equation}
Since $T$ is a contraction, $T^r$ is also a contraction for all $r\geq 1$. Hence $(I-T^{*r}T^r)\geq 0$ and $(I-T^r T^{*r})\geq 0$ for all $r\geq 1$.
This together with equations (\eqref{D_Tr} and \eqref{D_T*r} imply,
\begin{equation}\label{ker inq}
\begin{split}
\|F_r x\| &\leq \|P_n x\|, \ \ \text{and} \\
\|F_r^*x\| &\leq \|P_{n+kr}x\|
\end{split}
\end{equation}
holds for all $x\in \clh$ and all $r\geq 1$. Since $T$ is c.n.u, we will have (page 9, \cite{NF}),
\[
\left(\bigcap_{r\geq 1}\ker D_{T^r}\right)\bigcap\left(\bigcap_{r\geq 1}\ker D_{T^{*r}}\right)=\{0\}.
\]
Hence it follows by \eqref{ker D_T^r} and \eqref{ker D_T^*r}, that zero is the only common solution to the corresponding equalities in \eqref{ker inq}.
This settles the assertion $(4)$.

We will now prove the assertion $(3)$. By Douglas Lemma,
\begin{align*}
\begin{aligned}
\cld_T \subseteq \cld_{T^*}
\iff & D_T D_T^* \leq \lambda D_{T^*}D_{T^*}^* \text{ for some } \lambda\geq 0,
\\
\iff & D_T^2 \leq \lambda D_{T^*}^2
\\
\iff & I-T^*T \leq \lambda (I-TT^*)
\\
\iff & P_n- F_1^*F_1 \leq \lambda (P_{n+k}-F_1F_1^*), \text{ by \eqref{D_T=D_F_1} and \eqref{D_T*=D_F_1^*}}
\\
\iff &\lambda F_1F_1^*-F_1^*F_1 \leq \lambda P_{n+k}-P_n
\\
\iff & \lambda \|F_1^*x\|^2- \|F_1x\|^2 \leq \lambda \|P_{n+k}x\|^2-\|P_n x\|^2, \quad \forall x\in \clh.
\end{aligned}
\end{align*}
\end{proof}
Our next theorem shows that, any c.n.u contraction $T$ with finite indices and $\cld_T\subseteq\cld_{T^*}$ can be written as some finite rank perturbation of unilateral shift of a finite multiplicity. It is well known that a c.n.u contraction $T$ with $\cld_T=\{0\}$ is a pure isometry and hence is a unilateral shift of some multiplicity. Such operators are always analytic and hyponormal. To avoid this case, we will assume that\textsf{ the contraction $T$ in the below theorem is not an isometry.}
\begin{Theorem}\label{thm: char}
Let $T$ be a bounded linear operator on a Hilbert space $\clh$ and $n, k \in \N$. Then $T$ is a c.n.u contraction such that $\cld_T \subseteq\cld_{T^*}$ with $\dim\cld_T=n$, $\dim\cld_{T^*}<\infty$ and $\dim (\cld_{T^*}\ominus\cld_T)=k$, if and only if there exists an orthonormal basis$\{e_m\}_{m\geq 1}$ of $\clh$ with respect to which $T$ can be written as $T= S_k+F$, where $S_k$ is the unilateral shift of multiplicity $k$ and $F$ is a finite rank operator satisfying conditions $(1) - (4)$ in Proposition \ref{prop: char}.
\end{Theorem}
\begin{proof}
Let $T$ be a c.n.u contraction on $\clh$, with its defect spaces satisfying all the conditions as stated. It is well known (page 8, \cite{NF}) that $T: \clh\ominus\cld_T\longrightarrow \clh\ominus \cld_{T^*}$ is unitary. Since $T$ is c.n.u and $\dim \cld_T$ is a non zero finite number, $\cld_T$ must be a proper subspace of $\cld_{T^*}$. Hence $T: \clh\ominus \cld_T \longrightarrow \clh\ominus \cld_T$ is an isometry. Now,
\[
(\clh\ominus\cld_T)\ominus T(\clh\ominus\cld_T)=(\clh\ominus\cld_T)\ominus (\clh\ominus \cld_{T^*})=\cld_{T^*}\ominus \cld_T
\]
\NI and $\dim(\cld_{T^*}\ominus\cld_T)=k$. We show that $T$ on $\clh\ominus\cld_T$ is a pure isometry. On the contrary, assume $M\subseteq \clh\ominus\cld_T$ is such that $T|_{\clh\ominus\cld_T}$ on $M$ is unitary and $M$ reduces $T|_{\clh\ominus\cld_T}$. Note that $M=TM$ implies $M\subseteq T(\clh\ominus\cld_T)=\cld_{T^*}^{\perp}$. This shows that $\cld_{T^*}\subseteq M^{\perp}$. Since
\[
M^{\perp}=H\ominus M= \cld_T \oplus \Big((\clh\ominus\cld_T)\ominus M \Big),
\]
and $M$ reduces $T|_{\clh\ominus\cld_T}$, we have $T|_{\clh\ominus\cld_T}\Big((\clh\ominus\cld_T)\ominus M \Big)\subseteq (\clh\ominus\cld_T)\ominus M \subseteq M^{\perp}$. Again, $T(\cld_T)\subseteq \cld_{T^*}\subseteq M^{\perp}$. Hence  $T(M^{\perp})\subseteq M^{\perp}$. This shows that $M$ reduces the operator $T$ and $T|_M$ is unitary. Since $T$ is a c.n.u., $M$ must be zero. Hence, $T: \clh\ominus\cld_T\longrightarrow \clh\ominus\cld_T$ is a pure isometry, and by Wold decomposition, $T$ can be realized as a unilateral shift of multiplicity $k=\dim(\cld_{T^*}\ominus\cld_T)$. Let us denote $T|_{\clh\ominus \cld_T}$ by $S_k'$. Since $\dim\cld_T =n < \infty$ and $T$ maps $\cld_T$ into $\cld_{T^*}$, observe that $T|_{\cld_T}: \cld_T \longrightarrow \cld_{T^*}$ is a finite rank operator of rank at most $n$. Let us define $F_1: \clh \longrightarrow \clh$ by
\[
F_1 x=
\begin{cases}
Tx & \mbox{if } x\in\cld_T,
\\
0 & \mbox{if } x\in \cld_T^{\perp}.
\end{cases}
\]
Suppose, $\{e_1, e_2,\ldots, e_n\}$ and $\{e_{n+1}, e_{n+2},\ldots\}$ are the orthonormal bases of $\cld_T$ and $(\clh\ominus \cld_T)$ respectively, such that with respect to $\{e_{n+1}, e_{n+2},\ldots\}$, $T|_{\clh\ominus\cld_T}$ is $S_k'$. Clearly $\{e_m\}_{m\geq 1}$ is an orthonormal basis of $\clh$. Also, there exist scalars $x^i_j \in\C$ for $i=1, 2,\ldots n$, $j\geq 1$ such that, $F_1$ can be written as,
\[
F_1 (e_i)=
\begin{cases}
\sum_{j\geq 1}x^i_j e_j & \mbox{if } i\leq n,
\\
0 & \mbox{if } i\geq n+1.
\end{cases}
\]
The matrix representation of $T$ with respect to this orthonormal basis $\{e_m\}_{m\geq 1}$ is given by,
\begin{equation}\label{eqn: T matrix}
[T] = \begin{bmatrix}
x^1_1 & x^2_1 &\cdots &x^n_1 & 0 & 0 & 0 & \cdots
\\
x^1_2 & x^2_2 &\cdots & x^n_2 & 0 & 0 & 0 & \cdots
\\
\vdots & \vdots & \vdots & \vdots & \vdots &\vdots &\vdots & \vdots
\\
x^1_{n+k} & x^2_{n+k} &\cdots &x^n_{n+k} & 0 & 0 & 0 & \cdots
\\
x^1_{n+k+1} & x^2_{n+k+1} &\cdots &x^n_{n+k+1} & 1 & 0 & 0 & \cdots
\\
x^1_{n+k+2} & x^2_{n+k+2} &\cdots &x^n_{n+k+2} & 0 & 1 & 0 & \cdots
\\
x^1_{n+k+3} & x^2_{n+k+3} &\cdots &x^n_{n+k+3} & 0 & 0 & 1 & \cdots
\\
\vdots & \vdots & \vdots &\vdots &\vdots &\vdots & \vdots &\ddots
\end{bmatrix}.
\end{equation}
Since $T$ is a contraction, $T^*$ is also so. Hence, $\|T^*e_{n+k+l}\|\leq 1$ for all $l\geq 1$. This yields,
\[1+|x^1_{n+k+l}|^2+|x^2_{n+k+l}|^2+\cdots+|x^n_{n+k+l}|^2\leq 1.
\]
Hence $x^p_{n+k+l}=0$ for $p=1, 2, \ldots , n$ and $l\geq 1$. Now $F_1$ becomes,
\[
F_1 (e_i)=
\begin{cases}
\sum_{j= 1}^{n+k} x^i_j e_j & \mbox{if } i\leq n,
\\
0 & \mbox{if } i\geq n+1.
\end{cases}
\]
Now, the matrix representation of $T$ in \eqref{eqn: T matrix} (with respect to $\{e_m\}_{m\geq 1}$) reduces to
\begin{equation}\label{eqn: T new matrix}
[T] = \begin{bmatrix}
x^1_1 & x^2_1 &\cdots &x^n_1 & 0 & 0 & 0 & \cdots
\\
x^1_2 & x^2_2 &\cdots & x^n_2 & 0 & 0 & 0 & \cdots
\\
\vdots & \vdots & \vdots & \vdots & \vdots &\vdots &\vdots & \vdots
\\
x^1_k & x^2_k &\cdots & x^n_k & 0 & 0 & 0 & \cdots
\\
1+(x^1_{k+1}-1) & x^2_{k+1} &\cdots & x^n_{k+1} & 0 & 0 & 0 & \cdots
\\
\vdots & \vdots & \vdots &\vdots &\vdots &\vdots & \vdots &\vdots
\\
x^1_{n+k} & x^2_{n+k} &\cdots & 1+(x^n_{n+k}-1) & 0 & 0 & 0 & \cdots
\\
0 & 0 &\cdots & 0 & 1 & 0 & 0 & \cdots
\\
0 & 0 &\cdots & 0 & 0 & 1 & 0 & \cdots
\\
\vdots & \vdots & \vdots &\vdots &\vdots &\vdots & \vdots &\ddots
\end{bmatrix}.
\end{equation}
Let us define the operator $F: \clh \longrightarrow \clh$ by,
\[
F (e_i)=
\begin{cases}
x_1^ie_1+ x_2^i e_2+\cdots +(x_{i+k}^i-1)e_{i+k}+\cdots +x_{n+k}^i e_{n+k} & \mbox{if } i\leq n,
\\
0 & \mbox{if } i\geq n+1.
\end{cases}
\]
Then, with respect to $\{e_m\}_{m\geq 1}$, $T$ can be written as
\begin{equation}\label{T=S_k+F}
T=S_k+F
\end{equation}
where $S_k$ is the unilateral shift of multiplicity $k$ on $\clh$. Clearly, $F_1=F+S_kP_n$ and if we write $F_r$ for $r\geq 2$ as in \eqref{F_r}, then $T=S_k+F$ will satisfy all the assertions of Proposition \ref{prop: char}. Hence $S_k$ and $F$ will satisfy all the conditions $(1)-(4)$ as stated there.

Conversely, let there exists an orthonormal basis $\{e_m\}_{m\geq 1}$ on $\clh$ such that $T$ can be written as $T= S_k+F$, where $S_k$ acts on $\{e_m\}_{m\geq 1}$ as the unilateral shift of multiplicity $k$ and $F$ is a finite rank operator satisfying the conditions $(1)-(4)$ of the Proposition \ref{prop: char}. Then $S_k(e_m)=e_{m+k}$ for all $m\geq 1$ and there exist scalars $x^i_j\in\C, i=1,2,\ldots n$ and $ j= 1,2, \ldots n+k$ such that
\[
F (e_i)= (F_1-S_kP_n)(e_i)=
\begin{cases}
x_1^ie_1+ x_2^i e_2+\cdots +(x_{i+k}^i-1)e_{i+k}+\cdots +x_{n+k}^i e_{n+k} & \mbox{if } i\leq n,
\\
0 & \mbox{if } i\geq n+1,
\end{cases}
\]
\[
\text{where} \quad F_1 (e_i)=
\begin{cases}
x_1^ie_1+ x_2^i e_2+\cdots + x_{i+k}^i e_{i+k}+\cdots +x_{n+k}^i e_{n+k} & \mbox{if } i\leq n,
\\
0 & \mbox{if } i\geq n+1.
\end{cases}
\]
Now, with respect to $\{e_m\}_{m\geq 1}$, the matrix representation of $T$ will be,
\begin{equation}\label{eqn: T new matrix}
[T] = \begin{bmatrix}
x^1_1 & x^2_1 &\cdots &x^n_1 & 0 & 0 & 0 & \cdots
\\
x^1_2 & x^2_2 &\cdots & x^n_2 & 0 & 0 & 0 & \cdots
\\
\vdots & \vdots & \vdots & \vdots & \vdots &\vdots &\vdots & \vdots
\\
x^1_k & x^2_k &\cdots & x^n_k & 0 & 0 & 0 & \cdots
\\
x^1_{k+1} & x^2_{k+1} &\cdots & x^n_{k+1} & 0 & 0 & 0 & \cdots
\\
\vdots & \vdots & \vdots &\vdots &\vdots &\vdots & \vdots &\vdots
\\
x^1_{n+k} & x^2_{n+k} &\cdots & x^n_{n+k} & 0 & 0 & 0 & \cdots
\\
0 & 0 &\cdots & 0 & 1 & 0 & 0 & \cdots
\\
0 & 0 &\cdots & 0 & 0 & 1 & 0 & \cdots
\\
\vdots & \vdots & \vdots &\vdots &\vdots &\vdots & \vdots &\ddots
\end{bmatrix}.
\end{equation}
Proceeding exactly in the same way as in the proof of Proposition (\ref{prop: char}), it follows by condition $(4)$ that, $T$ is a c.n.u contraction. The condition $(1)$ implies $\dim\cld_T=n$, and then the inequality in $(3)$ says $\cld_T\subseteq\cld_{T^*}$. Lastly, the condition $\|F_r^* x\|\leq \|P_{n+kr}x\|$ with $r=1$ in $(4)$ shows $\dim\cld_{T^*}\leq n+k\leq\infty$. If possible, let $\dim\cld_{T^*}< n+k$. Then $\dim (\cld_{T^*}\ominus\cld_T)= k_1< k$ for some $k_1 \in\N$. Now as $T$ is a c.n.u contraction with $\cld_T\subseteq\cld_{T^*}<\infty$ such that $\dim\cld_T=n$ and $\dim (\cld_{T^*}\ominus\cld_T)= k_1$, the first part of the theorem immediately shows $T= S_{k_1}+F'$, for $S_{k_1}$ is the unilateral shift of multiplicity $k_1$ and $F'$ a finite rank operator with rank atmost $n$. Since the Fredholm index remains unchanged under finite rank perturbation, we must have $k_1=k$, which is a contradiction. Therefore, $\dim (\cld_{T^*}\ominus\cld_T)=k$.
\end{proof}
\begin{Remark}
A c.n.u contraction $T$ with $\cld_T\subseteq\cld_{T^*}$ need not be hyponormal nor analytic in general. The following example shows this.
\end{Remark}
\begin{Example}\label{n-hyp/analytic} Let $T$ be a bounded linear operator on $H^2(\D)$ defined by, $T(1)=T(z)=\frac{1}{2}$ and $T(z^m)= z^{m+1}$ for all $m\geq 2$.
The matrix representation of $T$ with respect to the orthonormal basis $\{z^m\}_{m\geq 1}$ is given by,

\[
[T] = \begin{bmatrix}
\frac{1}{2} & \frac{1}{2} & 0  & 0 & \cdots
\\
0 & 0 & 0 & 0 & \cdots
\\
0 & 0 & 0  & 0 & \cdots
\\
0 & 0 & 1 & 0 & \cdots
\\
0 & 0 & 0 & 1 & \cdots
\\
\vdots & \vdots & \vdots & \vdots & \ddots
\end{bmatrix}
\]
A straightforward computation will show that,

\[
[I-T^*T] = \begin{bmatrix}
\frac{3}{4} & -\frac{1}{4} & 0  & 0 & \cdots
\\
-\frac{1}{4} & \frac{3}{4} & 0 & 0 & \cdots
\\
0 & 0 & 0  & 0 & \cdots
\\
0 & 0 & 0 & 0 & \cdots
\\
\vdots & \vdots & \vdots & \vdots & \ddots
\end{bmatrix}.
\]
 Since the sub-matrix
$\begin{pmatrix}
\frac{3}{4} & -\frac{1}{4} \\
-\frac{1}{4} & \frac{3}{4}
\end{pmatrix}$  is positive definite and
$\det \begin{pmatrix}
\frac{3}{4} & -\frac{1}{4} \\
-\frac{1}{4} & \frac{3}{4}
\end{pmatrix}= \frac{1}{2}>0$, $T$ is a contraction with $\dim \cld_T= 2$ and $\cld_T= \overline{span}\{1, z\}$.

We next show that, $T$ satisfies the conditions $(1)-(4)$ of Proposition \ref{prop: char}. For this, we write $T$ (with resect to $\{z^m\}_{m\geq 0}$) as $T= F_1+S(I-P_2)$, where $S$ is the Hardy shift on $H^2(\D)$ and $F_1(1)=F_1(z)=\frac{1}{2}$ and $F_1(z^m)=0$ for $m\geq 2$. Clearly $F_1(I-P_2)=(I-P_3)F_1=0$ and hence the condition $(2)$ of the Proposition \ref{prop: char} is satisfied. By \eqref{F_r}, it follows that,

\[ F_r=
\begin{cases}
F_1  & \mbox{if } r=1, \\
F_1^r  & \mbox{if } r\geq 2.
\end{cases}
\]
Now for $x= \sum_{m=0}^{\infty}x_m z^m \in H^2(\D)$, we have $F_r(x)= \frac{1}{2^{r}}(x_0+x_1)$ $\forall$ $r\geq 1$. By Cauchy-Schwarz inequality,
$\frac{1}{2^{2r}}|x_0+x_1|^2 \leq \frac{1}{2^{2r-1}}(|x_0|^2+|x_1|^2)$ $\forall\ r\geq 1$. Hence $\|F_r x\|\leq \|P_2 x\|$ for all $x\in H^2(\D)$ and $r\geq 1$. Similarly, it can be shown, $\|F_r^* x\|\leq\|P_{2+r}x\|$ for all $x\in H^2(\D)$ and $r\geq 1$.

Let $x = \sum_{m=0}^{\infty}x_m z^m\in H^2(\D)$, be such that $\|F_r x\|=\|P_2 x\|$ and $\|F_r^* x\|=\|P_{2+r}x\|$ $\forall\ r\geq 1$. We show that $x=0$. First, $\|F_r x\|=\|P_2 x\|$ for all $r\geq 1$ implies, $\frac{1}{2^{2r}}|x_0+x_1|^2= |x_0|^2+|x_1|^2$ for all $r\geq 1$. Observe that, $x_0=x_1=0$ as $r\rightarrow \infty$. Therefore $x$ reduces to $x= \sum_{m=2}^{\infty}x_m z^m$.
Hence $\|F_r^* x\|=0$ for all $r \geq 1$. Now $\|F_r^* x\|=\|P_{2+r}x\|$ for all $r\geq 1$ implies,
\[
|x_2|^2+|x_3|^2+\cdots +|x_{1+r}|^2=0
\]
and hence $x_2=x_3=\cdots= x_{1+r}=0$, for all $r\geq 1$. This implies that $x=0$. Hence condition $(4)$ of the Proposition \ref{prop: char} holds.

Proceeding similarly as in the proof of Proposition \ref{prop: char}, one can show that, $I-F_1^*F_1= I-T^*T$. By the matrix form of the latter, $\text{rank  }(I-F_1^*F_1) =2$, and this is equivalent to satisfying condition $(1)$ of Proposition \ref{prop: char}.

Again, another simple computation will show that,

\[
[I-TT^*] = \begin{bmatrix}
\frac{1}{2} & 0 & 0  & 0 & \cdots
\\
0 & 1 & 0 & 0 & \cdots
\\
0 & 0 & 1  & 0 & \cdots
\\
0 & 0 & 0 & 0 & \cdots
\\
\vdots & \vdots & \vdots & \vdots & \ddots
\end{bmatrix},
\]

\NI and hence $\dim\cld_{T^*}=3$ and $\cld_{T^*}=\overline{span}\{1, z, z^2 \}$. Clearly $\cld_T\subseteq\cld_{T^*}$. Hence by Douglas lemma, there exists $\lambda\geq 0$ such that, $\lambda||F_1^*x||^2-||F_1x||^2\leq \lambda||P_3 x||^2-||P_2 x||^2$ holds for all $x\in H^2(\D)$ and this is condition $(3)$ in the Proposition \ref{prop: char}.
Hence, by the converse of Theorem~\ref{thm: char}, $T$ is a c.n.u., with $\cld_T\subseteq\cld_{T^*}$ such that $\dim\cld_T=2$ and $\cld_{T^*}\ominus\cld_T=1$.

We now show $T$ is not hyponormal. We begin by the computation (with respect to~$\{z^m\}_{m\geq 0}$),
\[
[T^*T-TT^*] = \begin{bmatrix}
-\frac{1}{4} & \frac{1}{4} & 0  & 0 & \cdots
\smallskip\\
\frac{1}{4} & \frac{1}{4} & 0 & 0 & \cdots
\\
0 & 0 & 1  & 0 & \cdots
\\
0 & 0 & 0 & 0 & \cdots
\\
\vdots & \vdots & \vdots & \vdots & \ddots.
\end{bmatrix}.
\]

Since $\la (T^*T-TT^*)(1), 1\ra= \la(-\frac{1}{4}+\frac{1}{4}z), 1\ra=-\frac{1}{4}< 0$, $(T^*T-TT^*)$ is not positive definite and hence is not hyponormal. Now $1 \in T^m H^2(\D)$ for all $m \geq 1$, implies that $T$ is not analytic.
\end{Example}

We will discuss the analyticity of a c.n.u. contraction in the setting of Theorem \ref{thm: char} in the next section; Now we discuss their hyponormality. Note for a hyponormal contraction $T$, $\cld_T \subseteq\cld_{T^*}$. But the converse is not true in general, see example \eqref{n-hyp/analytic} above. At the same time, condition $(3)$ in Theorem \ref{thm: char} suggests a condition for hyponormality. Indeed, we can show, under the assumption of Theorem \ref{thm: char}, $T$ is hyponormal if and only if $\lambda =1$.
\begin{Theorem}\label{Hyponormal}
Let T (and also $\cld_T$ and $\cld_{T^*}$) be as in the forward implication of Theorem~\ref{thm: char}. Then $T$ is hyponormal if and only if for all $x\in \clh$,
 $$\|F_1^*x\|^2-\|F_1x\|^2\leq \|P_{n+k}x\|^2-\|P_n x\|^2$$ holds, where $F_1$ is the finite rank operator associated with $T$ with respect to some orthonormal basis $\{e_m\}_{m\geq 1}$ as defined in \eqref{F_1'}.
\end{Theorem}
\begin{proof}
We know that $T$ is hyponormal if and only if $T^*T-TT^*\geq 0$. By Theorem \ref{thm: char}, $T$ can be written as $T=S_k+F$  with respect to some orthonormal basis $\{e_m\}_{m\geq 1}$ where $S_k$ and $F$ will satisfy the conditions $(1)-(4)$ of Proposition \ref{prop: char}.
Now,
\begin{align*}
\begin{aligned}
 T^*T-TT^* \geq 0 \iff & (I-TT^*)- (I-T^*T) \geq 0 \\
\iff & (P_{n+k}-F_1F_1^*)-(P_n-F_1^*F_1) \geq 0, (\text{ by \eqref{D_T=D_F_1}, \eqref{D_T*=D_F_1^*}})\\
\iff & F_1F_1^*- F_1^*F_1   \leq P_{n+k}-P_n \\
\iff & \|F_1^*x\|^2-\|F_1x\|^2 \leq \|P_{n+k}x\|^2-\|P_n x\|^2,\quad \forall x\in \clh.
\end{aligned}
\end{align*}
\end{proof}
We will now discuss a particular case of Theorems \ref{thm: char} in the setting of a c.n.u. contraction $T$ with $\dim\cld_T=1$. Unlike the other c.n.u. contractions, we will show that such $T$'s are always hyponormal.
\begin{Theorem}\label{1-Hypo}
Let $T$ be a bounded linear operator on a Hilbert space $\clh$. Then $T$ is a c.n.u contraction with $\cld_T\subseteq\cld_{T*}$ such that $\dim\cld_T=1$ and $\dim\cld_{T^*}=k+1>1$ if and only if there exists an orthonormal basis $\{e_n\}_{n\geq 0}$ on $\clh$, with respect to which $T$ can be written as $T=S_k+F$, where $S_k$ is the unilateral shift of multiplicity $k$ acting on $\{e_n\}_{n\geq 0}$, and $F$ has rank $\leq 1$ and is defined by
\begin{align*}
& F(e_0)=\alpha_0e_0+\alpha_1e_1+\cdots+(\alpha_k-1)e_k \text{ with }(\sum_{i=0}^{k}|\alpha_i|^2)<1 \text{ where } \{\alpha_j\}_{j=0}^k\in\C,\text{ and }\\ & F(e_n)=0 \ \forall \ n\geq 1.
\end{align*}
\end{Theorem}
\begin{proof}
Let $T$ be a c.n.u contraction on $\clh$ with $\cld_T\subseteq\cld_{T*}<\infty$ such that $\dim\cld_T=1$ and $\dim\cld_{T^*}=k+1$ for some $k> 0$. Theorem \ref{thm: char} yields an orthonormal basis $\{e_n\}_{n\geq 0}$ such that $T=S_k+F$, where $S_k(e_n)=e_{n+k}$ $\forall$ $n\geq 0$ and $F=F_1-S_kP_1$ with $F_1(e_0)=\alpha_0e_0+\alpha_1e_1+\cdots+\alpha_ke_k$ for some scalars $\alpha_0,\cdots,\alpha_k\in\C$ and $F_1(e_n)=0$ $\forall$ $n\geq 1$. Additionally, $F_1$ also satisfies the conditions $(1)-(4)$ of the Proposition \ref{prop: char}. Now, by condition $(4)$ of the same proposition with $x=e_0$ we get, $\|F_1 e_0\|\leq \|P_1e_0\|$ and this yields,
\begin{equation}\label{dim D_T=1}
|\alpha_0|^2+ |\alpha_1|^2 + \cdots + |\alpha_k|^2 \leq 1
\end{equation}
We will show that the strict inequality holds in \eqref{dim D_T=1}. Note that, the matrix representation of $T$ with respect to the orthonormal basis $\{e_n\}_{n\geq 0}$ is given by,
\begin{equation}\label{eqn: T (1) matrix}
[T] = \begin{bmatrix}
\alpha_0 &  0 & 0 &  \cdots
\\
\alpha_1 & 0 & 0 & \cdots
\\
\vdots & \vdots & \vdots & \vdots
\\
\alpha_k & 0 & 0 & \cdots
\\
0 & 1 & 0 & \cdots
\\
0 & 0 & 1 & \cdots
\\
\vdots & \vdots & \vdots &\ddots
\end{bmatrix}.
\end{equation}
By an easy computation, the matrix form of $I-T^*T$ with respect to the same basis is,
\begin{equation}\label{eqn: D_T (1) matrix}
[I-T^*T] = \begin{bmatrix}
1-(\sum_{i=0}^{k}|\alpha_i|^2) &  0 & 0 &  \dots
\\
0 & 0 & 0 & \dots
\\
0 & 0 & 0 & \dots
\\
\vdots & \vdots & \vdots &\ddots
\end{bmatrix}.
\end{equation}
Since $\dim \cld_T=1$, we must have $\sum_{i=0}^{k}|\alpha_i|^2< 1$.

For the converse part, let there exists an orthonormal basis $\{e_n\}_{n\geq 0}$, with respect to which, $T= S_k+F$ (for $F$ as in the statement) and can be represented as
\begin{equation}\label{eqn: T (2) matrix}
[T] = \begin{bmatrix}
\alpha_0 &  0 & 0 &  \cdots
\\
\alpha_1 & 0 & 0 & \cdots
\\
\vdots & \vdots & \vdots & \vdots
\\
\alpha_k & 0 & 0 & \cdots
\\
0 & 1 & 0 & \cdots
\\
0 & 0 & 1 & \cdots
\\
\vdots & \vdots & \vdots &\ddots
\end{bmatrix}.
\end{equation}
A similar computation to get the operator $(I-T^*T)$, as done in the first part, together with $\sum_{i=0}^{k}|\alpha_i|^2< 1$ would imply, $T$ is a contraction with $\dim\cld_T=1$. By \eqref{eqn: T (2) matrix}, we can write $T$ as $T= F_1+ S_k(I-P_1)$, where $F_1(e_0)= \alpha_0 e_0+\alpha_1 e_1+\cdots +\alpha_k e_k$ and $F_1(e_n)=0$ $\forall$ $n\geq 1$. We show that $F$ and $F_1$ satisfy conditions $(1)-(4)$ of the Proposition \ref{prop: char}. Clearly $F= F_1-S_kP_1$ and $F_1(I-P_1)=(I-P_{1+k})F_1=0$. Since $I-T^*T= P_1-F_1^*F_1$, we will have, $\text{rank}(P_1-F_1^*F_1)=1$. We now show that the inequality $(3)$ in Proposition \ref{prop: char} holds for $\lambda =1$.

Let $x=x_0e_0 + x_1e_1 + x_2e_2 + \cdots \in\clh$ be arbitrary. Then
\[
\|F_1^*x\|^2=\|(x_0\bar{\alpha_0}+x_1\bar{\alpha_1}+ \cdots + x_k\bar{\alpha_k})e_0\|^2= |x_0\bar{\alpha_0}+x_1\bar{\alpha_1}+ \cdots + x_k\bar{\alpha_k}|^2, \quad\text{and}
\]
\[
\|F_1x\|^2= \|x_0(\alpha_0e_0+\alpha_1e_1+\cdots + \alpha_k e_k)\|^2=|x_0|^2(|\alpha_0|^2+|\alpha_1|^2+ \cdots + |\alpha_k|^2).
\]
Therefore,
\begin{equation}\label{F_1*-F_1}
\|F_1^*x\|^2-\|F_1 x\|^2= |x_0\bar{\alpha_0}+x_1\bar{\alpha_1}+ \cdots + x_k\bar{\alpha_k}|^2- |x_0|^2 (|\alpha_0|^2+|\alpha_1|^2+ \cdots + |\alpha_k|^2).
\end{equation}
Now by Cauchy-Schwarz inequality,
\[
|x_0\bar{\alpha_0}+x_1\bar{\alpha_1}+ \cdots + x_k\bar{\alpha_k}|^2\leq (|x_0|^2+|x_1|^2+\cdots + |x_k|^2)(|\alpha_0|^2+|\alpha_1|^2+\cdots + |\alpha_k|^2).
\]
Hence by \eqref{F_1*-F_1}, and by $\sum_{i=0}^{k}|\alpha_i|^2 <1$,
\[
\|F_1^*x\|^2-\|F_1 x\|^2 \leq (|x_1|^2+|x_2|^2+\cdots+|x_k|^2)(|\alpha_0|^2+|\alpha_1|^2+\cdots+ |\alpha_k|^2)< |x_1|^2+|x_2|^2+\cdots+|x_k|^2.
\]
Next, since $|x_1|^2+|x_2|^2+\cdots+|x_k|^2= \|P_{1+k}x\|^2-\|P_1x\|^2$, we will have,
\begin{equation}\label{F_1hyp }
\|F_1^*x\|^2-\|F_1 x\|^2 \leq \|P_{1+k}x\|^2-\|P_1x\|^2\quad \forall\ x\in \clh.
\end{equation}
Setting
\[
F_r=
\begin{cases}
F_1, & \mbox{if } r=1, \\
F_1^r+S_k(I-P_1)F_1^{r-1}+\cdots+S_k^{r-1}(I-P_1)F_1 & \mbox{if } r\geq 2,
\end{cases}
\]
and proceeding exactly in the same way as in Proposition \ref{prop: char} we have, $I-T^{*r}T^r= P_1-F_r^*F_r$ $\forall$ $r\geq 1$ and $I-T^rT^{*r}=P_{1+kr}-F_rF_r^*$. As $T$ is a contraction, so is $T^r$ $\forall$ $r\geq 2$. Hence,
\begin{equation}\label{F_1}
\|F_rx\|\leq \|P_1 x\| \ \ \text{and}
\end{equation}
\begin{equation}\label{F_1*}
\|F_r^* x\|\leq \|P_{1+kr}x\|\ \ \text{for all } x \in \clh\ \text{ and }\ r\geq 1.
\end{equation}
We now show, the equalities in \eqref{F_1} and \eqref{F_1*} are attained only at $x= 0$ (for all $r$).

Let $x= x_0e_0+x_1e_1+\cdots + x_ne_n+\cdots \in \clh$ such that $\|F_rx\|=\|P_1x\|$ and $\|F_r^*x\|=\|P_{1+kr}x\|$ hold for all $r\geq 1$.
Now,
\begin{align*}
\begin{aligned}
 \|F_1x\| = \|P_1x\| \iff & \|x_0(\alpha_0e_0+\alpha_1e_1+\cdots+\alpha_ke_k)\|^2 = |x_0|^2 \\
\iff & |x_0|^2(|\alpha_0|^2+|\alpha_1|^2+\cdots+|\alpha_k|^2)=|x_0|^2.
\end{aligned}
\end{align*}
Since $(|\alpha_0|^2+|\alpha_1|^2+\cdots+|\alpha_k|^2)<1$, $x_0$ must be zero. Hence $x$ becomes $x= x_1e_1+x_2e_2+\cdots + x_ne_n+\cdots$.
Again,
\begin{align*}
\begin{aligned}
\|F_1^* x\| = \|P_{1+k}x\| \iff &\|(x_1\bar{\alpha_1}+x_2\bar{\alpha_2}+\cdots+x_k\bar{\alpha_k}) e_0\|^2 = \|x_1e_1+x_2e_2+\cdots + x_ke_k\|^2  \\
\iff & |x_1\bar{\alpha_1}+x_2\bar{\alpha_2}+\cdots+x_k\bar{\alpha_k}|^2 = |x_1|^2+|x_2|^2+\cdots+|x_k|^2.
\end{aligned}
\end{align*}
By Cauchy-Schwarz inequality,
\[
|x_1\bar{\alpha_1}+x_2\bar{\alpha_2}+\cdots+x_k\bar{\alpha_k}|^2 \leq (|x_1|^2+|x_2|^2+\cdots+|x_k|^2)(|\alpha_1|^2+|\alpha_2|^2+\cdots+|\alpha_k|^2).
\]
Since $\sum_{i=0}^{k}|\alpha_i|^2< 1$, we will have, $x_1=x_2=\cdots=x_k=0$. Therefore $x$ becomes,
$$x=x_{k+1}e_{k+1}+x_{k+2}e_{k+2}+\cdots+ x_n e_n+\cdots.$$
Now, for $r=2$, $F_2 = F_1^2+ S_k(I-P_1)F_1$. Now proceeding exactly in the similar way, the equation $\|F_2^*x\|=\|P_{1+2k}x\|$ will imply,
$x_{k+1}=x_{k+2}=\cdots= x_{2k}=0$.

Let $r\geq 2$. We show that, $x_0=x_1=\cdots= x_{rk}=0$ implies  $$x_{rk+1}=x_{rk+2}=\cdots=x_{(r+1)k}=0.$$
Let $x_0=x_1=\cdots= x_{rk}=0$. Then $x= x_{rk+1}e_{rk+1}+x_{rk+2}e_{rk+2}+\cdots$.
Recall, by \eqref{F_r},
\[
F_{r+1}^*= F_1^{*(r+1)}+F_1^{*r}(I-P_1)S_k^*+\cdots+ F_1^{*2}(I-P_1)S_k^{*(r-1)}+F_1^*(I-P_1)S_k^{*r} \ \forall\ r\geq 2.
\]
For $1\leq p\leq r-1$,
\begin{equation}\label{F_1*(I-P_1)}
\begin{split}
F_1^*(I-P_1)S_k^{*p}(x)&= F_1^*(I-P_1)S_k^{*p}(x_{rk+1}e_{rk+1}+x_{rk+2}e_{rk+2}+\cdots)\\
&= F_1^*(I-P_1)(x_{rk+1}e_{rk+1-pk}+x_{rk+2}e_{rk+2-pk}+\cdots).
\end{split}
\end{equation}
\begin{align*}
\begin{aligned}
\text{Note},\ \  p  \leq r-1\implies  pk  \leq rk-k \quad (\text{as }k\geq 1)
\implies & k  \leq rk-pk\\
\implies & rk-pk+l  \geq k+l\geq k+1 \quad \forall\ l \geq 1.
\end{aligned}
\end{align*}
So, $F_1^*(e_{rk+l-pk})=0$ for all $l \geq 1$ and hence by (\ref{F_1*(I-P_1)}),
\[
F_1^*(I-P_1)(x_{rk+1}e_{rk+1-pk}+x_{rk+2}e_{rk+2-pk}+\cdots)=0.
\]
Consequently, $F_1^{*(r-p+1)}(I-P_1)S_k^{*p}x=0$ for $1\leq p\leq (r-1)$ and thereby,
\[
F_{r+1}^*x= F_1^*(I-P_1)S_k^{*r}(x).
\]
\[
\begin{split}
\text{So,}\ \  F_{r+1}^*x &= F_1^*(I-P_1)S_k^{*r}(x_{rk+1}e_{rk+1}+ x_{rk+2}e_{rk+2}+ \cdot\cdot+ x_{(r+1)k}e_{rk+k}+x_{(r+1)k+1}e_{(r+1)k+1}+\cdot\cdot)
\\
 &= F_1^*(I-P_1)(x_{rk+1}e_1+ x_{rk+2}e_2+ \cdots+ x_{(r+1)k}e_k+ x_{(r+1)k+1}e_{k+1}+\cdots)
\\
&= (x_{rk+1}\bar{\alpha_1}+ x_{rk+2}\bar{\alpha_2}+ \cdots+ x_{(r+1)k}\bar{\alpha_k})e_0,
\end{split}
\]
and $\|F_{r+1}^*x\|=\|P_{(r+1)k+1}x\|$ together imply,
\[
|x_{rk+1}\bar{\alpha_1}+x_{rk+2}\bar{\alpha_2}+\cdots+ x_{(r+1)k}\bar{\alpha_k}|^2 = |x_{rk+1}|^2+|x_{rk+2}|^2+\cdots+ |x_{(r+1)k}|^2.
\]

\NI Again, Cauchy-Schwarz inequality and $\sum_{i=0}^{k}|\alpha_i|^2<1$ will imply, $x_{rk+1}=x_{rk+2}=\ldots= x_{(r+1)k}=0$. Now by the induction hypothesis, $x_m=0$ for all $m\geq 0$ and so $x=0$. This shows that the condition $(4)$ of Proposition \ref{prop: char} is satisfied.
Therefore, by the converse part of Theorem \ref{thm: char}, $T$ is a c.n.u contraction with $\cld_T\subseteq\cld_{T^*}<\infty$ such that $\dim\cld_T=1$ and $\dim\cld_{T^*}=k+1$.
\end{proof}
An immediate consequence of Theorems \ref{Hyponormal} and \eqref{1-Hypo}, and inequality \eqref{F_1hyp } is:
\begin{Corollary}\label{D_T=1}
If $T$ is a c.n.u contraction on a Hilbert space $\clh$ with finite dimensional defect spaces such that $\cld_T \subseteq\cld_{T^*}$ and $\dim\cld_T=1$, then $T$ is hyponormal.
\end{Corollary}

Note, when $\alpha_0=0$ in Theorem \ref{1-Hypo}, $T$ is analytic. Indeed, $\alpha_0=0$ implies, $T\clh \subseteq \overline{span}\{e_1, e_2, e_3, \ldots\}$ so that for $n\geq 1$, $T^{n+1}\clh \subseteq \{e_{nk+1}, e_{nk+2},\ldots\}$ and hence $\bigcap_{m\geq 1}T^m \clh=\{0\}$. However if $\alpha_0 \neq 0$, $T$ need not be analytic. For example, if $\alpha_1=\alpha_2=\cdots=\alpha_k=0$, then $\alpha_0$ is an eigenvalue of $T$ corresponding to eigenvector $e_0$. Hence $\{e_0\} \in \bigcap_{m\geq 1} T^m \clh$.
In the next section, we will explore when such operators can be analytic in a more general setting.

\newsection{Analyticity}\label{sec: Analyticity}

In this section, we will give the criteria for the analyticity of c.n.u contractions $T$ on $\clh$ with finite indices such that $\cld_T\subseteq\cld_{T^*}$. We begin by fixing $\dim\cld_T=n+1$ and $\dim(\cld_{T^*}\ominus\cld_T)=k$, for some $n\geq 0$ and $k\geq 1$. By Theorem \ref{thm: char}, an orthonormal basis $\{e_m\}_{m\geq 0}$ exists on $\clh$ with respect to which the matrix representation of $T$ is given by,

\begin{equation}\label{eqn: T Analytic matrix}
[T] = \begin{bmatrix}
a_{00}& a_{01} & a_{02} &\cdots & a_{0n}& 0 & 0 & \cdots
\\
a_{10}& a_{11} & a_{12} &\cdots & a_{1n} & 0 & 0 & \cdots
\\
\vdots & \vdots & \vdots & \vdots & \vdots &\vdots &\vdots & \vdots
\\
a_{n0}& a_{n1} & a_{n2} &\cdots & a_{nn} & 0 & 0 & \cdots
\\
a_{n+1,0}& a_{n+1,1} & a_{n+1,2} &\cdots & a_{n+1,n} & 0 & 0 & \cdots
\\
\vdots & \vdots & \vdots & \vdots & \vdots &\vdots &\vdots & \vdots
\\
a_{n+k,0}& a_{n+k,1} & a_{n+k,2} &\cdots & a_{n+k,n} & 0 & 0 & \cdots
\\
0 & 0 & 0 & \cdots & 0 & 1 & 0 & \cdots
\\
0 & 0 & 0 & \cdots & 0 & 0 & 1 & \cdots
\\
\vdots & \vdots & \vdots &\vdots &\vdots &\vdots & \vdots &\ddots
\end{bmatrix}
\end{equation}
for some suitable scalars $a_{ij}$, $i= 0,1,2,\ldots, n+k$ and $j= 0,1,2,\ldots, n$.

Via the canonical unitary map $U: \clh\longrightarrow H^2(\D)$, defined by $U(e_m)=z^m, m\geq 0$, $T$ can be viewed as an operator on the Hardy space $H^2(\D)$ having the same matrix representation \eqref{eqn: T Analytic matrix} with respect to $\{z^m\}_{m\geq 0}$.

Let us consider the finite rank operators $A$ and $B$ on $H^2(\D)$ whose matrix representations with respect to $\{z^m\}_{m\geq 0}$, respectively are,
\begin{equation}\label{eqn: A Analytic matrix}
[A] = \begin{bmatrix}
a_{00}& a_{01} & a_{02} &\cdots & a_{0n}& 0 &  \cdots
\\
a_{10}& a_{11} & a_{12} &\cdots & a_{1n} & 0 & \cdots
\\
\vdots & \vdots & \vdots & \vdots & \vdots &\vdots & \vdots
\\
a_{n0}& a_{n1} & a_{n2} &\cdots & a_{nn} & 0 & \cdots
\\
0 & 0 & 0 &\cdots & 0 & 0 & \cdots
\\
\vdots & \vdots & \vdots &\vdots &\vdots &\vdots & \ddots
\end{bmatrix},\ \
[B] = \begin{bmatrix}
0 & 0 & 0 &\cdots & 0 & 0 & \cdots
\\
0 & 0 & 0 &\cdots & 0 & 0 & \cdots
\\
\vdots & \vdots & \vdots & \vdots & \vdots &\vdots & \vdots
\\
0 & 0 & 0 &\cdots & 0 & 0 & \cdots
\\
a_{n+1,0}& a_{n+1,1} & a_{n+1,2} &\cdots & a_{n+1,n} & 0 & \cdots
\\
a_{n+2,0}& a_{n+2,1} & a_{n+2,2} &\cdots & a_{n+2,n} & 0 &\cdots
\\
\vdots & \vdots & \vdots & \vdots & \vdots &\vdots & \vdots
\\
a_{n+k,0}& a_{n+k,1} & a_{n+k,2} &\cdots & a_{n+k,n} & 0 & \cdots
\\
0 & 0 & 0 & \cdots & 0 & 0 & \cdots
\\
0 & 0 & 0 & \cdots & 0 & 0 & \cdots
\\
\vdots & \vdots & \vdots &\vdots &\vdots &\vdots & \ddots
\end{bmatrix}.
\end{equation}
 Then $T$ can be written as $$T=A+B+S_k(I-P_{n+1}),$$ where $S_k$ is the unilateral shift of multiplicity $k$ on $H^2(\D)$
and $P_{n+1}$ is the orthogonal projection onto the space $\overline{span}\{1, z, z^2,\ldots, z^n\}$.
We will first show via inducting on $r\geq 2$, that
\begin{equation}\label{T^r=A+B}
T^r= A^r+BA^{r-1}+S_k B A^{r-2}+S_k^2 B A^{r-3}+\cdots+ S_k^{r-2}B A+ S_k^{r-1}B + S_k^r(I-P_{n+1}).
\end{equation}
Assume that $r=2$. Then,
\begin{equation*}
\begin{split}
 T^2 =&\ [A+B+S_k(I-P_{n+1})][A+B+S_k(I-P_{n+1})]\\
=&\ A^2+AB+AS_k(I-P_{n+1})+BA+B^2+BS_k(I-P_{n+1})+S_k(I-P_{n+1})A\\
&+S_k(I-P_{n+1})B+S_k(I-P_{n+1})S_k(I-P_{n+1}).
\end{split}
\end{equation*}
\medskip
Since the operators $AB,  AS_k(I-P_{n+1}),  B^2,  BS_k(I-P_{n+1}),  S_k(I-P_{n+1})A$ are all zero and since $S_k(I-P_{n+1})S_k(I-P_{n+1})=S_k^2(I-P_{n+1})$ and $S_k(I-P_{n+1})B=S_k B$, we have
\[
T^2= A^2 + BA + S_k B+S_k^2(I-P_{n+1}).
\]
Hence, the equality in \eqref{T^r=A+B} holds for $r=2$.
\smallskip

For the induction step, let the equality hold for some $r\geq 3$. Then,
\[
\begin{split}
T^{r+1}= T^r T
=&\ [A^r+BA^{r-1}+S_k B A^{r-2}+S_k^2 B A^{r-3}+\cdots+ S_k^{r-2}B A + S_k^{r-1}B\\
&\ + S_k^r(I-P_{n+1})]\ [A+B+S_k(I-P_{n+1})] \\
=&\ A^{r+1}+BA^r+S_k BA^{r-1}+S_k^2BA^{r-2}+\cdots+ S_k^{r-2}BA^2+S_k^{r-1}BA\\
&\   +S_k^r(I-P_{n+1})B+S_k^{r+1}(I-P_{n+1})  \\
=&\  A^{r+1}+BA^r+S_k BA^{r-1}+ S_k^2BA^{r-2}+\cdots + S_k^{r-1}BA+S_k^r B+ S_k^{r+1}(I-P_{n+1}),
\end{split}
\]
and hence by the induction hypothesis, the proof of \eqref{T^r=A+B} follows.
Using this equality, we investigate the analyticity as well as the non-vanishing eigenvalues of $T$, by exploring several cases of $A$ and $B$ (appearing in \eqref{T^r=A+B} for $T$).

\textsf{We first discuss the case when $A$ is nilpotent.} Then $A^{n+1}=0$ and hence $T^{n+1}H^2(\D)\subseteq \{e_{n+1}, e_{n+2},\ldots\}$, which further implies,
\[
T^{n+l+1}H^2(\D)\subseteq \{e_{n+kl+1}, e_{n+kl+2},\ldots\}\quad \forall\  l\geq 1.
\]
\[
\text{So,}\quad \bigcap_{l\geq 1}T^{n+l+1}H^2(\D)\subseteq \bigcap_{l\geq 1}\{e_{n+kl+1}, e_{n+kl+2},\ldots\}=\{0\},
\]
and hence $\bigcap_{m\geq 1}T^m H^2(\D)=\{0\}$. Therefore, \textsf{if $A$ is nilpotent, $T$ is analytic.}

\NI Note that, with respect to $\{z^m\}_{m\geq 0}$, $A$ can be written as, $ A=\begin{pmatrix}A_1 & 0  \\0 & 0\end{pmatrix}$ where $A_1(z^m)=\sum_{i=0}^{n}a_{im}z^i$, for $0\leq m\leq n$. Clearly,

\begin{enumerate}
\item $A$ is nilpotent if and only if $A_1$ is nilpotent.
\item $A$ has a non-zero eigenvalue if and only if $A_1$ has a non-zero eigenvalue.
\end{enumerate}

 Now \textsf{we assume that $A$ is not nilpotent.} This is equivalent to saying that $A_1$ is not nilpotent. We discuss the (sub-) cases $B=0$ and $B\neq 0$, separately.

\textsf{Suppose $B=0$.} Then there exists $\lambda(\neq 0)\in\C$ and $h=(h_i)_{i=0}^n (\neq 0)\in\C^{n+1}$ such that $A_1h=\lambda h$. Now if we consider the sequence $h'=(h_0, h_1, \ldots, h_n, 0,0,\ldots)$, then it is easy to see that $Th'=\lambda h'$. This will imply, $h'\in T^m H^2(\D)$ for all $m\geq 1$, and hence \textsf{$T$ is not analytic.}

\textsf{Now, we assume that $B\neq 0$ (and $A$ is not-nilpotent).} Let $\lambda (\neq 0)\in\C$ be an eigenvalue of $T$. Then there exists a non-zero vector $h=\sum_{m=0}^{\infty}h_m z^m \in H^2(\D)$ such that $Th=\lambda h$. So,
\newpage
\[
\begin{bmatrix}
a_{00}& a_{01} & a_{02} &\cdots & a_{0n}& 0 & 0 & \cdots
\\
a_{10}& a_{11} & a_{12} &\cdots & a_{1n} & 0 & 0 & \cdots
\\
\vdots & \vdots & \vdots & \vdots & \vdots &\vdots &\vdots & \vdots
\\
a_{n0}& a_{n1} & a_{n2} &\cdots & a_{nn} & 0 & 0 & \cdots
\\
a_{n+1,0}& a_{n+1,1} & a_{n+1,2} &\cdots & a_{n+1,n} & 0 & 0 & \cdots
\\
\vdots & \vdots & \vdots & \vdots & \vdots &\vdots &\vdots & \vdots
\\
a_{n+k,0}& a_{n+k,1} & a_{n+k,2} &\cdots & a_{n+k,n} & 0 & 0 & \cdots
\\
0 & 0 & 0 & \cdots & 0 & 1 & 0 & \cdots
\\
0 & 0 & 0 & \cdots & 0 & 0 & 1 & \cdots
\\
\vdots & \vdots & \vdots &\vdots &\vdots &\vdots & \vdots &\ddots
\end{bmatrix}
\begin{bmatrix}
h_0 \\
h_1\\
\vdots\\
h_n \\
h_{n+1} \\
h_{n+2}\\
\vdots\\
h_{n+k}\\
h_{n+k+1} \\
h_{n+k+2}\\
\vdots
\end{bmatrix}
=
\lambda
\begin{bmatrix}
h_0 \\
h_1\\
\vdots\\
h_n \\
h_{n+1} \\
h_{n+2}\\
\vdots\\
h_{n+k}\\
h_{n+k+1} \\
h_{n+k+2}\\
\vdots
\end{bmatrix},\quad\text{and hence}
\]
\begin{equation}\label{T eigenvalue}
\begin{bmatrix}
a_{00}h_0+a_{01}h_1+a_{02}h_2+\cdots+ a_{0n}h_n \\
a_{10}h_0+a_{11}h_1+a_{12}h_2+\cdots+ a_{1n}h_n\\
\vdots\\
a_{n0}h_0+a_{n1}h_1+a_{n2}h_2+\cdots+ a_{nn}h_n \\
a_{n+1,0}h_0+a_{n+1,1}h_1+a_{n+1,2}h_2+\cdots+ a_{n+1,n}h_n \\
a_{n+2,0}h_0+a_{n+2,1}h_1+a_{n+2,2}h_2+\cdots+ a_{n+2,n}h_n \\
\vdots\\
a_{n+k,0}h_0+a_{n+k,1}h_1+a_{n+k,2}h_2+\cdots+ a_{n+k,n}h_n \\
h_{n+1} \\
h_{n+2}\\
\vdots\\
h_{n+k}\\
h_{n+k+1}\\
h_{n+k+2}\\
\vdots\\
h_{n+2k}\\
\vdots
\end{bmatrix}
=\lambda
\begin{bmatrix}
h_0 \\
h_1\\
\vdots\\
h_n \\
h_{n+1} \\
h_{n+2}\\
\vdots\\
h_{n+k}\\
h_{n+k+1} \\
h_{n+k+2}\\
\vdots\\
h_{n+2k}\\
h_{n+2k+1}\\
h_{n+2k+2}\\
\vdots\\
h_{n+3k}\\
\vdots
\end{bmatrix}.
\end{equation}
From the above expression, it follows via successive iterations,
\begin{equation}\label{eqn: h iterates}
h_{n+pk+l}=\frac{h_{n+l}}{\lambda^p} \text{ for }p\geq 0,\quad l= 1,\ldots k.
\end{equation}
Now suppose $h_m=0$ for all $m=0,1,\ldots, n$. By \eqref{T eigenvalue}, $\lambda \neq 0$ forces $h_{n+l}=0$ for all $l= 1, \ldots, k$ and by \eqref{eqn: h iterates}, $h_{n+pk+l}=0$ for all $l= 1,\ldots, k$ and $p\geq 1$. So $h=0$, contradicting the assumption that $h$ is an eigenvector.
Therefore, $h_m\neq 0$ for at least one $m\in \{0,1, \ldots,n\}$.

From equation \eqref{T eigenvalue}, we have $A_1\begin{bmatrix}
h_0 \\
\vdots\\
h_n
\end{bmatrix}
=
\lambda \begin{bmatrix}
h_0 \\
\vdots\\
h_n
\end{bmatrix}$, i.e., $\lambda\neq 0$ is also an eigenvalue of $A_1$.

We now show that,
\begin{equation}\label{B eigen vector}
a_{n+l,0}h_0+a_{n+l,1}h_1+a_{n+l,2}h_2+\cdots +a_{n+l,n}h_n=0, \text{ for all } l= 1,2,\ldots k.
\end{equation}

If possible, let there exist $l\in \{1,2,\ldots, k\}$, such that the L.H.S of \eqref{B eigen vector} is non-zero, and for the sake of definiteness we assume $l=1$. By \eqref{T eigenvalue} this implies $h_{n+1}\neq 0$, and by \eqref{eqn: h iterates} $h_{n+pk+1}\neq 0$ for all $p\geq 0$. Since $h\in H^2(\D)$, the sequence $h_{n+1}(1,\frac{1}{\lambda}, \frac{1}{\lambda^2}, \ldots)$ is square summable. But this is a contradiction as $0<|\lambda|\leq \|T\|\leq 1$.

Conversely, if an eigenvector $(h_0, h_1,\cdots, h_n)$ corresponding to an eigenvalue $\lambda\neq 0$ of $A_1$ satisfies $\sum_{i=0}^{n}a_{n+l,i}h_i=0$ for all $l= 1, \ldots , k$, then \eqref{T eigenvalue} says $Th=\lambda h$ for $h=(h_0, h_1, h_2,\ldots, h_n, 0, 0, \ldots)$. So, $\lambda \neq 0$ is also an eigenvalue of $T$.
Formally, this summarizes:
\begin{Lemma}\label{B-eigen vector}
Let $T$ be a c.n.u contraction on a Hilbert space $\clh$ such that $\cld_T\subseteq\cld_{T^*}$ with $\dim\cld_{T^*}<\infty$. Let $\dim\cld_T=n+1$ and $\dim(\cld_{T^*}\ominus\cld_T)=k$ for some $n\geq 0$ and  $k\geq 1$. Assume that, $T$ can be represented as \eqref{eqn: T Analytic matrix} with respect to some orthonormal basis $\{e_n\}_{n\geq 0}$ on $\clh$. Then, $T$ has no non-zero eigenvalue if and only if for every eigenvector $(h_0, h_1, \ldots, h_n)$ corresponding to any non-zero eigenvalue of the sub-matrix
\[
A_1=\begin{bmatrix}
a_{00}& a_{01} & a_{02} &\cdots & a_{0n}
\\
a_{10}& a_{11} & a_{12} &\cdots & a_{1n}
\\
\vdots & \vdots & \vdots & \vdots & \vdots
\\
a_{n0}& a_{n1} & a_{n2} &\cdots & a_{nn}
\end{bmatrix},
\]
there exists at least one $l\in \{1, 2, \ldots, k\}$, such that, $(\sum_{r=0}^{n}a_{n+l,r}h_r)\neq 0$.
\end{Lemma}
From the beginning of this section, we dealt with the case of $A$ being a  contraction. Now suppose only the submatrix $A_1$ of $A$ is a contraction. Then arguing similarly we have:
\begin{Corollary}\label{Non-contractive}
Let $T$ be a bounded linear operator on a Hilbert space $\clh$ with matrix representation \eqref{eqn: T Analytic matrix} of which only the submatrix $[A_1]$ of $[A]$ in \eqref{eqn: A Analytic matrix} is a contraction. Then the Lemma \ref{B-eigen vector} will hold as good for any such operator $T$.
\end{Corollary}
We will now come to the main theorem of this section.
\begin{Theorem}\label{Analytic}
Let $T$  be a bounded linear operator on $H^2(\D)$ with matrix representation \eqref{eqn: T Analytic matrix} with respect to the basis $\{z^m\}_{m\geq 0}$. Let the submatrix
$A_1$ of $T$ (as in the above lemma) be a contraction. T is analytic if and only if $T$ has no non-zero eigenvalue.
\end{Theorem}
\begin{proof}
Let $T$ be analytic. If possible let $\lambda (\neq 0)\in\C$ be an eigenvalue of $T$. Then there exists $h\neq 0\in H^2(\D)$ such that $T^m h=\lambda^m h$ $\forall$ $m\geq 1$ and $h\in \bigcap_{m\geq 1}T^m H^2(\D)$, a contradiction.

Conversely, let $T$ have no non-zero eigenvalue; indeed such situations are possible in the light of Lemma \ref{B-eigen vector}. Now, two cases can arise:
\begin{enumerate}
\item $A_1$ is nilpotent.
\item $A_1$ has non-zero eigenvalues.
\end{enumerate}

\NI \textbf{Case 1:} Let $A_1$ be nilpotent. Then, as we have seen earlier, $T$ is analytic.\\
\NI \textbf{Case 2:} Let $A_1$ has non-zero eigenvalues.
\\
Since $A_1 \in M_{n+1}(\C)$ is a finite matrix, by Schur decomposition and the standard theory of matrix (\cite{HJ}, \cite{HK}), $A_1$ is similar to a lower triangular matrix with its eigenvalue $\lambda_0, \lambda_1,\ldots, \lambda_n$ on the main diagonal. Let there exist scalars $b_{ij}, b_{ij}', c_{ij}$ for $i,j = 0,1, 2,\ldots n$ such that
\[W=
\begin{bmatrix}
b_{00}& b_{01} & b_{02} &\cdots & b_{0n}
\\
b_{10}& b_{11} & b_{12} &\cdots & b_{1n}
\\
\vdots & \vdots & \vdots & \vdots & \vdots
\\
b_{n0}& b_{n1} & b_{n2} &\cdots & b_{nn}
\end{bmatrix}\in GL_{n+1}(\C),
\quad
W^{-1}=
\begin{bmatrix}
b_{00}' & b_{01}' & b_{02}' &\cdots & b_{0n}'
\\
b_{10}' & b_{11}' & b_{12}' &\cdots & b_{1n}'
\\
\vdots & \vdots & \vdots & \vdots & \vdots
\\
b_{n0}' & b_{n1}' & b_{n2}' &\cdots & b_{nn}'
\end{bmatrix},\ \ \text{and}
\]
\[W^{-1}A_1W=
\begin{bmatrix}
\lambda_0 & 0 & 0 &\cdots & 0
\\
c_{10}& \lambda_1 & 0 &\cdots & 0
\\
c_{20} & c_{21}& \lambda_2 &\cdots & 0
\\
\vdots & \vdots & \vdots & \vdots & \vdots
\\
c_{n0}& c_{n1} & c_{n2} &\cdots & \lambda_n
\end{bmatrix}.
\]

Without loss of generality, we may assume that the non-zero eigenvalues correspond to the first columns of $W^{-1}A_1W$, i.e., if $\lambda_0, \lambda_1,\ldots, \lambda_q$ with $(q\leq n)$ are the only non-zero eigenvalues of $A_1$, then $\lambda_0, \lambda_1,\ldots, \lambda_q $ appear on the first $(q+1)$ columns of $W^{-1}A_1W$. Clearly the operator with matrix representation (with respect to $\{z^m\}_{m\geq 0}$)
$\begin{pmatrix}
W & 0 \\
0 & I
\end{pmatrix}(=R)$ on $H^2(\D)$ is invertible with inverse
$\begin{pmatrix}
W^{-1} & 0 \\
0 & I
\end{pmatrix}(=R^{-1})$. Let $T_1= R^{-1}TR$. Then $T_1$ is a bounded linear operator on $H^2(\D)$. Since $T_1^m = R^{-1}T^m R$ for all $m\geq 1$, $T$ is analytic if and only if $T_1$ is analytic. We will show that $T_1$ is analytic. There exists suitable scalars $c_{n+l, j}$ for $l=1,2,\ldots, k$ and $j= 0,1, 2,\ldots, n$ such that with respect to $\{z^m\}_{m\geq 0}$, $T_1$ can be represented as,
\[
\begin{bmatrix}
\lambda_0 & 0 & 0 &\cdots & 0 & 0 & 0 & \cdots
\\
c_{10}& \lambda_1 & 0 &\cdots & 0 & 0 & 0 & \cdots
\\
c_{20}& c_{21} & \lambda_2 &\cdots & 0 & 0 & 0 & \cdots
\\
\vdots & \vdots & \vdots & \vdots & \vdots &\vdots &\vdots & \vdots
\\
c_{n0}& c_{n1} & c_{n2} &\cdots & \lambda_n & 0 & 0 & \cdots
\\
c_{n+1,0}& c_{n+1,1} & c_{n+1,2} &\cdots & c_{n+1,n} & 0 & 0 & \cdots
\\
c_{n+2,0}& c_{n+2,1} & c_{n+2,2} &\cdots & c_{n+2,n} & 0 & 0 & \cdots
\\
\vdots & \vdots & \vdots & \vdots & \vdots &\vdots &\vdots & \vdots
\\
c_{n+k,0}& c_{n+k,1} & c_{n+k,2} &\cdots & c_{n+k,n} & 0 & 0 & \cdots
\\
0 & 0 & 0 & \cdots & 0 & 1 & 0 & \cdots
\\
0 & 0 & 0 & \cdots & 0 & 0 & 1 & \cdots
\\
\vdots & \vdots & \vdots &\vdots &\vdots &\vdots & \vdots &\ddots
\end{bmatrix}.
\]
First, we assume $0\leq q<n$.
Clearly, $\overline{T_1H^2(\D)}=\overline{span} \{p_0, p_1, \ldots, p_n ;  z^{n+k+1}, z^{n+k+2}, \ldots\}$~where,
\[
p_r(z)=
\begin{cases}
\lambda_rz^r+c_{r+1,r}z^{r+1}+\cdot\cdot+c_{n,r}z^n+c_{n+1,r}z^{n+1}+c_{n+2,r}z^{n+2}+\cdot\cdot+ c_{n+k, r}z^{n+k}, & \mbox{if } 0\leq r< n \\
\lambda_nz^n+c_{n+1,n}z^{n+1}+c_{n+2,n}z^{n+2}+\cdots+c_{n+k,n}z^{n+k}, & \mbox{if }r=n.
\end{cases}
\]
Infact, $\overline{T_1\{1, z, z^2, \ldots, z^n\}}$ is spanned by $\{p_0, p_1, p_2, \ldots, p_n\}$. Note that $p_0, p_1, p_2, \ldots, p_q$ are linearly independent, and the non-zero elements (if any) of the set $\{p_{q+1}, p_{q+2}, \ldots, p_n\}$ are polynomials with $z^{ q+2}$ as a factor.
For any $m\geq 1$, $p_{r,m}(z)=T_1^mp_r(z) = \lambda_r^{m+1} z^rp_r'(z)$ for all r with $0\leq r\leq q$ and some polynomials $p_r'$ such that $p_r'(0)=1$. Also, $p_{r,m}(z)$ are linearly independent for $0\leq r\leq q$ and $m\geq 1$.
Let $m_1$ be the smallest positive integer such that non-zero elements (if any) of the set $ X:=\{T_1^{m_1}(p_{q+l}(z)): 1\leq l\leq (n-q)\}$ are polynomials with $z^{n+2}$ as a factor. Suppose, there are some non-zero elements in the set $X$. Without loss of generality, we may assume that all the elements of $X$ are non-zero.

Let us choose $m> m_1$. Then there exists $r_m \in\N$ such that $m= r_m+m_1$; clearly $r_m$ increases with $m$.
The polynomials,
$p_{q+l,m}(z)=T_1^m(p_{q+l}(z))$, for $1\leq l\leq (n-q)$, have a factor $z^{n+2+kr_m}$.

Note that, for $1\leq l\leq n-q $, $\deg p_{q+l,m}(z)\leq n+(m+1)k$.
Now, for any $m\geq m_1$,
\[
\overline{T_1^{m+1}\{1, z, z^2,\ldots, z^n\}}=\overline{span}\{p_{r,m}(z) (0\leq r\leq q);\quad p_{q+l,m}(z)(1\leq l\leq n-q)\},
\]
where, $p_{r,m}(z)$ are linearly independent for $(0\leq r\leq q)$. We can find elements from the set $\{p_{q+l,m}(z)(1\leq l\leq n-q)\}$ to form a basis for $\overline{T_1^{m+1}\{1, z, z^2,\ldots, z^n\}}$ containing the set $\{p_{r,m}(z), (0\leq r\leq q)\}$. For simplicity (as there will be no harm), let us assume that the set $\{p_{r,m}(z) (0\leq r\leq q);\quad p_{q+l,m}(z)(1\leq l\leq n-q)\}$ itself is linearly independent and form a basis for $\overline{T_1^{m+1}\{1, z, z^2,\ldots, z^n\}}$. Note that, this basis is not necessarily orthonormal.

Hence for $m\geq m_1$, the below set forms a basis for $\overline{T_1^{m+1}H^2(\D)}$:
\[
\{p_{r,m}(z) : 0\leq r\leq q\}\bigcup\{ p_{q+l,m}(z) : 1\leq l\leq n-q\}\bigcup \{z^{n+(m+1)k+1}, z^{n+(m+1)k+2}, \ldots\}.
\]

If possible let $h(\neq 0)\in H^2(\D)$, such that $h \in \bigcap_{p\geq 1}\overline{T_1^pH^2(\D)}$. Then $[h]_{T_1} \subseteq \bigcap_{p\geq 1}\overline{T_1^pH^2(\D)}$. Two subcases can arise: For $0\leq i\leq q$,

\begin{enumerate}[(a)]
\item $\la h, z^i\ra=0$
\item $\la h, z^i\ra \neq 0$.
\end{enumerate}

\NI \textbf{Subcase (a):} $\la h, z^i\ra=0$  for $0\leq i\leq q$.

We can write, $h= z^{q+l_1}(t_0+t_1z+t_2z^2+\cdots)$ for some $l_1\geq 1$ and $\sum_{m_1=0}^{\infty}t_{m_1}z^{m_1}$ with $t_0\neq 0$.
Corresponding to this $(q+l_1)$, we can choose $m(> m_1)$ large enough, so that $ q+l_1< n+2+k r_m$. For this $m$, $h\in \overline{T_1^{m+1}H^2(\D)}$.

Hence there exist scalars $\alpha_r (0\leq r\leq q)$, $\beta_l (1\leq l\leq n-q)$; and $\{\gamma_{l_2}, l_2\geq 1\}$ such that
\begin{equation}\label{h}
\begin{split}
h =&\  (\alpha_0p_{0,m}+\alpha_1p_{1,m}+\alpha_2p_{2,m}+\cdots+\alpha_q p_{q,m})+(\beta_1p_{q+1,m}+\beta_2p_{q+2,m}+\cdots \\
&\ \cdots+\beta_{n-q}p_{n,m})
  + \sum_{l_2\geq 1}\gamma_{l_2}z^{n+(m+1)k+l_2}.
\end{split}
\end{equation}
As $\Big(q+l_1< n+2+kr_m< n+(m+1)k+l_2\Big)$ and $\la h, z^{q+l_1}\ra=t_0\neq 0$, by equation (\ref{h}), $\alpha_i\neq 0$ for some $i\in \{0,\ldots,q\}$. Again, $q+l_1\geq q+1$ yields, $\alpha_0\lambda_0^{m+1}=0$. Since $\lambda_0$ is non-zero, we must have $\alpha_0=0$. Similarly $\alpha_r\lambda_r^{m+1} =0$ for $1\leq r\leq q$ and will imply $\alpha_r=0$ for $1\leq r\leq q$. This is a contradiction. Hence $\bigcap_{p\geq 1}\overline{T_1^pH^2(\D)}=\{0\}$ and consequently $\bigcap_{p\geq 1}T_1^pH^2(\D)=\{0\}$.

\NI \textbf{Subcase (b):} $\la h,z^i\ra\neq 0$ for some $i\in \{0,1,2,\ldots, q\}$.

We may assume that $i$ is the smallest such integer. We can write $h= z^i(t_0+t_1 z+t_2 z^2+\cdots)$ for some $\sum_{r=0}^{\infty}t_rz^r\in H^2(\D)$ such that $t_0\neq 0$.
Since $\lambda_i$ is not an eigenvalue of $T_1$, $h_1= T_1h-\lambda_i h$ is non-zero and $h_1 \in [h]_{T_1}\subseteq\bigcap_{p\geq 1}\overline{T_1^p H^2(\D)}$.

We can write $h_1= z^{i+j}(t_0'+t_1'z+t_2'z^2+\cdots)$ for some $j\geq 1$, where $t_0'\neq 0$.

\smallskip

If $i+j\geq (q+1)$ then $\la h_1,z^i\ra=0$ $\forall$ $i=0,1,\ldots,q$, and we will proceed as in subcase (a).

If $i+j\leq q$, we will consider the non-zero element $T_1h_1- \lambda_{i+j}h_1$ in $[h]_{T_1}$. As $q$ is finite, only after finitely many steps, we will find some $h' (\neq 0)\in\bigcap_{p\geq 1}\overline{T_1^p H^2(\D)}$ such that $\la h', z^i\ra=0$ for all $i=0,1,\ldots q$. Again, we will proceed as in subcase (a) and conclude that $T_1$ is analytic.

\textsf{If there exists $m_1'\in \mathbb{N}$ such that the  elements in the set $X'=\{T_1^{m_1'}p_{q+l}(z): 1\leq l\leq n-q\}$ are all zero, the subcases (a) and (b) will follow similarly and we will get the same conclusion.}

In the remaining case $q=n$, the polynomials $p_{r,m}(z)=T_1^m p_r(z)$, for $0\leq r\leq n$ and $m\geq 1$ are linearly independent and form a basis of $\overline{T_1^{m+1}\{1, z, z^2,\ldots, z^n\}}$. Hence the set $\{p_{r, m}(z),(0\leq r\leq n); z^{n+(m+1)k+l},l\geq 1\}$ is linearly independent and forms a basis of $T_1^{m+1}H^2(\D)$. Now proceeding similarly as before, the analyticity of $T_1$ will follow.
\end{proof}
\begin{Corollary}
Let $T$ on $\clh$ be a non-isometric c.n.u. contraction with $\cld_T\subseteq\cld_{T^*}$ and $\dim\cld_{T^*}<\infty$. Then $T$ is analytic if and only if it has no non-zero eigenvalue.
\end{Corollary}
\begin{proof}
Follows by Theorem \ref{thm: char} and Theorem \ref{Analytic}.
\end{proof}
\begin{Corollary}\label{1-analytic}
Let $T$ be a bounded linear operator on $H^2(\D)$ such that $T=S_k+F$ where~ $S_k$ acts on $\{z^n\}_{n\geq 0}$ as the unilateral shift of multiplicity $1\leq k <\infty$ and $F$ is defined by,
\[
F(z^n)=\begin{cases}
\alpha_0+\alpha_1 z+\cdots+(\alpha_k-1)z^k & \mbox{if } n=0, \\
0 & \mbox{if } n\geq 1,
\end{cases}
\]
for $\{\alpha_j\}_{j=0}^k\in\C$ and $0<|\alpha_0|\leq 1$. $T$ is analytic if and only if $\alpha_j\neq 0$ for some $1\leq j\leq k$.
\end{Corollary}
\begin{proof}
Follows by Lemma \ref{B-eigen vector}, Corollary \ref{Non-contractive}, and Theorem \ref{Analytic}.
\end{proof}
\begin{Corollary}
Let $T$ be a c.n.u. contraction with $\dim\cld_T=1$, $\cld_T\subseteq\cld_{T^*}$ and $\dim\cld_{T^*}<\infty$ (as in Theorem \ref{1-Hypo}). $T$ is analytic, if and only if either $\alpha_0=0$ or $\alpha_j\neq 0$ for some $1\leq j\leq k$.
\end{Corollary}
\begin{proof}
The proof is a direct consequence of the Corollary \ref{1-analytic} and the discussion before Lemma~\ref{B-eigen vector}, as $\alpha_0=0$ can be considered as a nilpotent matrix of order one over $\C$.
\end{proof}
\vspace{0.1in}

\noindent\textbf{Acknowledgement:}
The author is thankful to Prof. E. K. Narayanan for helpful discussions, suggestions, and corrections throughout the work.
The research of the author is supported by the IoE-IISc fellowship at the Indian Institute of Science, Bangalore, India.

\bibliographystyle{amsplain}

\end{document}